\documentclass{test}

\usepackage{amsmath,amsfonts,amssymb,color,enumerate}

\def\red{black} 

\def\C{\mathbb C}  
\def\D{\mathbb D}  
\newcommand{\IdObject}{\mathsf{E}}
\def\HH{{\bf H}}  
\def\SS{{\bf S}} 
\def\I{\mathbb I}  
  
\def\K{\mathbb K}  
\def\L{\mathfrak L}
\def\N{\mathbb N}  
\def\P{\mathbb P}  
\def\R{\mathbb R}  
\def\RR{\mathcal R}  
\def\Ss{\mathbb S}  
  
\def\Z{\mathbb Z}  
 
\def\phi{\varphi}
\def\epsilon{\varepsilon}
\def\.{{\cdot}}  
\def\<{\langle}  
\def\>{\rangle}  
\def\({\big(}   
\def\){\big)}      
\def\til{\widetilde}
\def\hat{\widehat}

\def\<{\langle} 
\def\>{\rangle}
\def\defi{\stackrel{\rm def}{=}}

\def\ssk{\smallskip}  
\def\msk{\medskip}   
\def\bsk{\bigskip}   
\def\nin{\noindent}  
\def\cen{\centerline} 
\font\bfone=cmbx10 scaled\magstep1

\DeclareMathOperator{\Hom}{Hom}
\DeclareMathOperator{\Der}{Der}
\DeclareMathOperator{\Span}{span}
\DeclareMathOperator{\colim}{colim}

\DeclareMathOperator{\Spec}{Spec}

\DeclareMathOperator{\id}{id}

\DeclareMathOperator{\ev}{ev}
\DeclareMathOperator{\im}{im}
\DeclareMathOperator{\quot}{quot}
\DeclareMathOperator{\inc}{inc} 
\DeclareMathOperator{\inv}{inv} 
\DeclareMathOperator{\conv}{conv} 
\DeclareMathOperator{\supp}{supp} 

\def\lead{\leaders\hbox to 1.5ex{\hss${.}$\hss}\hfill}
\def\arr{\hbox to 60pt{\rightarrowfill}}
\def\larr{\hbox to 60pt{\leftarrowfill}}
\def\mapdown#1{\Big\downarrow\rlap{$\vcenter{\hbox{$\scriptstyle#1$}}$}}
\def\lmapdown#1{\llap{$\vcenter{\hbox{$\scriptstyle#1$}}$}\Big\downarrow}
\def\mapright#1{\smash{\mathop{\arr}\limits^{#1}}}
\def\lmapright#1{\smash{\mathop{\arr}\limits_{#1}}}
\def\mapleft#1{\smash{\mathop{\larr}\limits^{#1}}}
\def\lmapleft#1{\smash{\mathop{\larr}\limits_{#1}}}
\def\mapup#1{\Big\uparrow\rlap{$\vcenter{\hbox{$\scriptstyle#1$}}$}}
\def\lmapup#1{\llap{$\vcenter{\hbox{$\scriptstyle#1$}}$}\Big\uparrow}

\def\tilbar #1{\overline{\til{\hbox{$#1$}}}} 
\def\bartil #1{\til{\bar{\hbox{$#1$}}}}

\def\mdvertspace{\vphantom{A_{A_A}}}
\def\muvertspace{\vphantom{A^{A^A}}}
\def\mdddvertspace{\vphantom{A_{A_{A_{A_A}}}}}
\def\muuuvertspace{\vphantom{A^{A^{\int^{o}a}}}}
\theoremstyle{defn}
\theoremstyle{thm}

\long\def\alert#1{\parindent2em\smallskip\hbox to\hsize%
{\hskip\parindent\vrule%
\vbox{\advance\hsize-2\parindent\hrule\smallskip\parindent.4\parindent%
\narrower\noindent#1\smallskip\hrule}\vrule\hfill}\smallskip\parindent0pt}

\title{The Pro-Lie Group Aspect of  Weakly Complete Algebras\\
           and Weakly Complete Group Hopf Algebras}
                  
\author{Rafael Dahmen and Karl Heinrich Hofmann}

\lastname{Dahmen and Hofmann}

\msc{22e15, 22e65, 22e99}  


\address{%
Rafael Dahmen\\	
Karlsruher Institut f\"ur\\ 
Technologie (KIT)	\\
76131 Karlsruhe, Germany\\
rafael.dahmen@kit.edu}

\address{%
Karl Heinrich Hofmann\\       
Fachbereich Mathematik\\      
Technische Universit\"at Darmstadt\\
Schlossgartenstra{\ss}e 7\\
64289 Darmstadt, Germany\\       
\hbox{hofmann@mathematik.tu-darmstadt.de}}

\begin{document}

\maketitle

\vglue-20pt

\cen{\it Dedicated to Joachim Hilgert on the Occasion of his 60th birthday}

\bsk
\begin{abstract} A weakly complete vector space over $\K=\R$ or $\K=\C$
is isomorphic to $\K^X$ for some set $X$ algebraically and topologically.
The significance of this type of topological vector spaces is illustrated
by the fact that the underlying  vector space of the Lie algebra of any 
pro-Lie group is weakly complete. 
In this study,  weakly complete real or complex associative algebras are studied
because they are necessarily
projective limits of finite dimensional algebras. The group of units $A^{-1}$
of a weakly complete algebra $A$ is a pro-Lie group with the associated
topological Lie algebra $A_{\rm Lie}$ of $A$ as Lie algebra and the
globally defined exponential function $\exp\colon A\to A^{-1}$ as
the exponential function of $A^{-1}$. 
With each topological group, a weakly complete group
algebra $\K[G]$ is associated functorially so that the functor
$G\mapsto \K[G]$ is left adjoint to $A\mapsto A^{-1}$. The group
algebra $\K[G]$ is a weakly complete Hopf algebra. If $G$ is compact,
the $\R[G]$  contains $G$ as the set of grouplike elements. 
The category of all real Hopf algebras $A$ with a compact group
of grouplike elements whose linear span is dense in $A$ is shown
to be equivalent to the category of compact groups. The group algebra
$A=\R[G]$ of a compact group $G$ contains a copy of the Lie algebra $\L(G)$
in $A_{\rm Lie}$; it also contains  a copy of the Radon 
measure algebra $M(G,\R)$. The
dual of the group algebra $\R[G]$ is the Hopf algebra ${\mathcal R}(G,\R)$
of representative functions of $G$. The rather straightforward duality 
between vector spaces
and weakly complete vector spaces thus becomes the basis 
of a duality ${\mathcal R}(G,\R)\leftrightarrow \R[G]$ and thus
 yields a new aspect of Tannaka duality.  
\end{abstract}

\bigskip

\cen{\bfone Introduction}

\msk

\nin
        A weakly complete vector space is a topological vector space
        isomorphic in the category of topological vector space{\color{\red}s}
        to ${\color{\red}\R}^X$ for some set $X$ with the Tychonoff topology.
	These vector spaces appear prominently in any Lie theory of 
        pro-Lie groups,  for instance, of  all compact groups, or even of
        all locally compact groups which are 
        compact modulo connectivity.
        This text will present a general study of
	weakly complete {\color{\red}associative} unital {\it algebras} 
        over the fields of real and complex 
	numbers, including eventually all weakly complete Hopf algebras. 
        It will turn out that all weakly complete algebras  are automatically
        projective limits of finite dimensional ones, and their 
        groups of units (that is, multiplicatively invertible elements)
        are always pro-Lie groups. This theory permits the introduction
        and study of weakly complete topological group algebras 
        over the real and the
        complex field  and their major properties. In particular, they
        are automatically topological  Hopf algebra{\color{\red}s}.
        Regarding compact groups, for example, we shall see
        that the weakly complete group algebra $\R[G]$ of a compact group $G$
        not only contains $G$ but also {\color{\red}its} Lie algebra 
        $\L(G)$ and {\color{\red}its}
        (Radon) measure algebra $M(G,\R)$.  {\color{\red} Indeeed let $\mathcal G$
        denote the category of compact groups and $\mathcal H$
        the category of weakly complete
        real Hopf algebras $A$ whose set $G$ of grouplike elements is 
        compact such that its linear span $\Span(G)$ is dense in $A$. 
        One of our central results is that $\mathcal G$ and $\mathcal H$ are 
        equivalent categories.  The dualtity between the category of
        abstract real vector spaces and the category of weakly complete
        vector spaces is of the Pontryagin type, is natural, and is
        easily accessible as
        has been detailed in the literature. Thus it is natural and relatively
        easy to pass from the category $\mathcal H$ to its dual category
        of real (abstract) Hopf algebras. In this process, the dual of $\R[G]$
        emerges as the  real (abstract) 
        Hopf algebra $\RR(G,\R)$ of reprsentative functions. 
        While this so-called {\it Tannaka duality} 
        between a compact group $G$ and the function vector space $\RR(G,\R)$
        is standard in the literature but also although somewhat arduous
        to follow through, 
        the natural equivalence $G\leftrightarrow \R[G]$
        between the categories $\mathcal G$ and $\mathcal H$ 
        however, now appears as a  natural link that newly enlightens 
        Tannaka duality.}

\bsk

      \section {Preliminaries: Vector Spaces and their Duality}
In this section we collect the preliminary concepts which secure
the language we use {\color{\red}for dealing} with the structure 
of a technically simple
concept, namely{\color{\red},} that of topological {\color{\red}associative}
unital algebras whose
underlying vector space is what we {\color{\red}shall} call weakly complete. 
A crucial
fact therefore is the duality theory of weakly complete vector spaces.

\bsk

\subsection{Duality Reviewed} For the sake of completeness 
we recall the formal definition
of a duality. We assume sufficient familiarity with category theory.
 A convenient reference is \cite{book3}, pp.\ 745--818.
\msk
 Let
 $\cal A$ and $\cal B$ be categories and consider two
contravariant functors (see \cite{book3}, Definition 3.25) 
$F\colon{\cal A}\to{\cal B}$ and $U\colon{\cal B}\to {\cal A}$.

\begin{Definition}\quad
(a)\quad  It is  said that $F$ 
and $U$ are {\it adjoint on the right} if for each 
object $A$ in $\cal A$ there is a natural $\cal A$-morphism 
$\eta_A\colon A\to UFA$ such that for each $\cal A$-morphism 
$f\colon A\to UB$ there is a 
unique $\cal B$-morphism $f'\colon B\to FA$ such that 
$f=(Uf')\circ \eta_A$.

 $$
\begin{matrix}&{\cal A}&&\hbox to 7mm{} &{\cal B}\\ 
\noalign{\vskip3pt}
\noalign{\hrule}\cr
\noalign{\vskip3pt}%
  \mdvertspace A&\mapright{\eta_A}&UFA&\hbox to 7mm{} &FA\\
\lmapdown{\forall f}&&\mapdown{Uf'}&\hbox to 7mm{}&\mapup{\exists! f'}\\
\muvertspace UB&\lmapright{\id}&UB&\hbox to 7mm{}&B.\\
\end{matrix}
 $$

\bsk

The natural morphism $\eta_A$ is called the {\it evaluation morphism}.

\nin (b) {\it The categories $\cal A$ and $\cal B$ are said to be dual}
if there are two contravariant functors $F$ and $U$, adjoint on the
right, such that
$\eta_A$ is an isomorphism for all objects $A$ of $\cal A$. 
\end{Definition}

\msk

\nin In the general literature, authors frequently 
 speak of ``duality'' as soon as they
have two ``naturally defined'' contravariant functors 
that are adjoint on the right without
the postulate that the natural morphism $\eta_A$ is \
an isomorphism in all cases---for example if they consider
the category of Banach spaces as dual to itself under the passage 
of a Banach space to the dual Banach space. 

\bsk

       \subsection{Vector Spaces Reviewed}
It is a basic fact that a perfect 
duality theory exists between the category of
real or complex vector spaces $\cal V$ and the linear 
maps between them on the one hand
and the so-called weakly complete 
topological vector spaces 
$\cal W$ and the continuous linear maps between them
on the other
(see \cite{dahmen, probook, book3, hofmor}).
Yet this fact is less commonly mentioned than one might
expect. 
In the following let $\K$ denote the field 
of either the real or the complex numbers. 
Most of the duality works with arbitrary topological fields, 
but for applications in the theory of, say, compact groups, 
a restriction to these fields suffices.

\begin{Definition} \label{d:wvect} A topological 
$\K$-vector space $V$ is called {\it weakly 
complete} if there is a vector space 
$E$ such that $V$ is isomorphic as a topological 
vector space to the space $E^*\defi\Hom(E,\K)$
of linear scalar valued functionals endowed 
with the topology of pointwise convergence
from the containment $\Hom(E,\K)\subseteq \K^E$.
\end{Definition}

\nin Together with all continuous linear maps between 
them, weakly complete topological vector spaces form
a category denoted $\cal W$.

In particular, since every $\K$-vector space 
$E$ is isomorphic to the direct sum 
$\K^{(J)}$ for some set $J$
due to the existence of a basis via the
Axiom of Choice we have from these definitions
the following insight:

\begin{Remark} \quad Every weakly complete vector space is
 isomorphic as a topological vector space to $\K^J$
for some set $J$.  
\end{Remark}

\msk 

\nin According to this remark the cardinality of $J$ is the only 
characteristic invariant of a weakly complete vector space.

\msk
 
Conversely, if a weakly complete vector space $V$ is given,
a vector space $E$ such as it is required by 
Definition \ref{d:wvect},
is easily obtained by considering the (algebraic) vector 
space of its topological dual
$V'=\Hom(V,\K)$ of all continuous linear functionals.
Indeed, we have:

\msk

\begin{Theorem} \label{th:vect} 
The categories $\cal V$ and $\cal W$ are dual to
each other with the contravariant functors 
$$E\mapsto E^*:{\cal V}\to {\cal W},\hbox{ respectively, }
V\mapsto V':{\cal W}\to {\cal V}.$$
\end{Theorem}

\nin For a vector space $E$, the evaluation morphism
$\eta_E\colon E\to (E^*)'$, for $v\in E$ and 
$\omega\in (E^*)'$, is given
by $\eta_E(v)(\omega)=\omega(v)$. 
In an analogous way one obtains the
evaluation morphism $\eta_V\colon V\to (V')^*$.

\msk

For proofs see \cite{book3}, Theorem 7.30, and \cite{probook}, 
Appendix 2: weakly complete topological vector spaces.

This duality theorem is one of the ``good'' ones for the category
$\cal W$ of weakly complete topological vector spaces since it
can be used for the understanding of the category $\cal W$
by transporting statements in and about it to well known 
{\color{\red}purely algebraic} statements
on vector spaces and linear maps.

\bsk
        \subsection{Monoidal categories}
The category $\cal V$ has a ``multiplication,'' namely, the tensor
product $(E,F)\mapsto E\otimes F$. Indeed this product
is associative and commutative and it has any one dimensional
vector space as identity object (represented conveniently by
 $\K$ itself). In order to meet the demands of general
category theory, there is a formalism around these concepts
which allows us to speak of {\it commutative monoidal
categories} as is exemplified by a discussion in \cite{book3},
Appendix 3, Paragraphs A3.61---A3.92, or in \cite{hofbig}, Section 2.
{\it Associativity} of the product $\otimes$ is implemented
by a natural isomorphism 
$\alpha_{EFG}\colon E\otimes(F\otimes G)\to (E\otimes F)\otimes G$
and {\it commutativity} by a natural isomorphism 
$\kappa_{EF}\colon E\otimes F\to F\otimes E$. (Cf.\ \cite{book3},
pp.~787ff.)
The simplest example of a commutative monoidal category is the
category of sets and functions with the cartesian binary product
$(X,Y)\mapsto X\times Y$ with the singleton sets as identity
elements.

\msk

If a functor between two monoidal categories respects the available
monoidal structures, it will be called {\it multiplicative}.
(See \cite{book3}, Definition A3.66.)

\msk

For the present duality theory of vector spaces it is relevant
that not only the category $\cal V$ of vector spaces
has a tensor product, but that the category $\cal W$ of weakly
complete vector spaces has a tensor product 
$(V,W)\mapsto V\otimes W$ as well (see \cite{heil} (1997) and
\cite{dahmen} (2007)). The tensor product of weakly complete
vector spaces has universal properties which are completely 
analogous to those well known for vector spaces. Indeed,
we have the following proposition:
\msk
\begin{Proposition} \label{p:multi} The respective tensor products endow each of
$\cal V$ and $\cal W$ with the structure of a commutative monoidal
category, and the contravariant functors 
$$E\mapsto E^*\colon {\cal V}\to{\cal W}\hbox{ and }
 V\mapsto V' \colon {\cal W}\to{\cal V}$$
are multiplicative. 
\end{Proposition}

\msk (See e.g.\ \cite{dahmen}). 
In particular, there are natural isomorphisms 
$$(E_1\otimes E_2)^*\cong (E_1)^*\otimes(E_2)^*\hbox{ and } 
(V_1\otimes V_2)' \cong (V_1)'\otimes(V_2)'.$$

     \subsection{Monoids and Comonoids---Algebras and Coalgebras}
Let us consider a commutative monoidal category $({\cal A}, \otimes)$.
Good examples for the purpose of the discussion are the categories
of sets with the cartesian product or, alternatively, the category 
of compact spaces and the (topological) cartesian product.

A {\it multiplication} on an object $A$ of $\cal A$ is an 
$\cal A$-morphism $m\colon A\otimes A\to A$. Its 
{\it associativity} can be expressed
in terms of $\cal A$-morphisms in a commutative diagram

\def\arr{\hbox to 30pt{\rightarrowfill}}
\def\larr{\hbox to 30pt{\leftarrowfill}}
\def\ot{\otimes}

$$\begin{matrix}%
A\ot(A\ot A)&\mapright{\alpha_{AAA}}&(A\ot A)\ot A\\
\lmapdown{\id_A\ot m}&&\mapdown{m\otimes \id_A}\\
A\otimes A && A\otimes A\\
\lmapdown{m}&&\mapdown{m}\\
A&\lmapright{\id_a}&A.\\
\end{matrix}
$$

\msk\nin
The multiplication $m$ is called 
{\it commutative} if the following diagram commutes:

$$
\begin{matrix}%
\mdvertspace A\ot A&\mapright{\kappa_{AA}}&A\ot A\\
\lmapdown{m}&&\mapdown{m}\\
\muvertspace A&\lmapright{\id_A}&\hphantom{.}A.\\
\end{matrix}
$$

\msk

(Cf.\ \cite{book3}, Definition A3.62).

A set $X$ with an associative multiplication is commonly called
a semigroup, a topological space $C$ with an associative
multiplication is a topological semigroup. 

A category
$\cal A$ is said to have identity objects $\IdObject$ (such as singletons in the case
of the categories of sets or topological spaces), if for each object $A$
there are natural isomorphisms $\iota_A\colon \IdObject\ot A\to A$ and 
$\iota'_A\colon A\ot \IdObject\to A$. 

If an object $A$ has a multiplication $m\colon A\ot A\to A$, then 
an identity of $(A,m)$ is a morphism $u\colon \IdObject\to A$ 
such that the following diagram commutes:

$$
\begin{matrix}
\mdvertspace \IdObject\ot A&\mapright{u\ot\id_A}&A\ot A&%
         \mapleft{\id_A\ot u}&A\ot \IdObject\\
\lmapdown{\iota_A}&&\lmapdown{m}&&\mapdown{\iota'_A}\\
\muvertspace A&\lmapright{\id_A}&A&\lmapleft{\id_A}&\hphantom{.}A.\\
\end{matrix}
$$

\msk
\begin{Definition} \label{d:monoid} 
In a commutative monoidal category 
$\cal A$, an object $A$ with an associative multiplication
$m\colon A\otimes A\to A$ and an identity $u\colon \IdObject\to A$
is called a {\it monoid} in $\cal A$. 
\end{Definition}

\msk

If we take for $(\cal A, \otimes)$ the category $(\cal V, \otimes)$ 
of vector spaces with the tensor product then 

\nin
{\it a monoid in $(\cal V,\otimes)$ is precisely a unital associative
algebra over $\K$,}  

\noindent
(in other words, an associative $\K$-algebra with identity). 
Usually, the multiplication $m\colon A\otimes A\to A$ is 
expressed in the form 
$$(a,b)\mapsto ab\defi m(a\otimes b):A\times A\to A,$$
while the identity $u\colon \R\to A$ gives rise to the 
element $1_A=u(1)$ satisfying $a1_A=1_Aa=a$ for all $a\in A$.
 
In exactly the same spirit we consider
the category $\cal W$ of weakly complete vector spaces. Here

\nin
{\it a monoid in $\cal W$ is precisely a unital associative
weakly complete topological algebra over $\K$.} 

    \section{Weakly Complete Associative Unital Algebras}

It is perhaps useful for our purpose to emphasize the utter simplicity
of this concept in a separate definition:

\begin{Definition} \label{d:basic-definition}
  A {\it weakly complete unital algebra} is an
associative algebra $A$ over $\K$ with identity, whose underlying
vector space is weakly complete, and whose multiplication
$(a,b)\mapsto ab:A\times A\to A$ is continuous.
\end{Definition}

\msk
\begin{Example} \label{e:3.2} (1a)\quad The product algebra $\K^J$ for 
a set $J$  
(with componentwise operations) is a weakly complete
unital algebra.

\ssk

(1b) More generally, if $\{A_j:j\in J\}$ is any family of finite 
dimensional unital $\K$-algebras, then $A=\prod_{j\in J} A_j$
is a weakly complete unital algebra over $\K$.

\ssk

(1c) Even more generally: Let $J$ be a directed set and 
$$\{p_{jk}\colon A_k\to A_j\quad|\quad j\le k, \hbox{ in }J\}$$
a projective system of morphisms of finite dimensional unital $\K$-algebras.
Then the projective limit $A=\lim_{j\in J}A_j$ is a weakly complete
unital $\K$-algebra. (See \cite{book3}, Definition 1.25ff., pp.17ff.,
also \cite{probook}, pp.77ff.) 

\msk
(2)\quad The algebra $\K[[X]]$ of all formal power series
$a_0+a_1X+a_2X^2+\cdots$ with $(a_0,a_1,a_2,\dots)\in\K^{\N_0}$ and the
topology making 
$$(a_n)_{n\in\N_0}\mapsto \sum_{n=0}^\infty a_nX^n:\K^{\N_0}\to \K[[X]]$$
an isomorphism of weakly complete vector spaces is a weakly complete
unital algebra. 

\msk

(3)\quad Let $V$ be any weakly complete topological vector space.
Endow $A=V\times\K$ with componentwise addition and the multiplication
$$ (v_1,r_1)(v_2,r_2)=(r_1\.v_2+r_2\.v_1,r_1r_2).$$
Then $A$ is a weakly complete unital algebra. 
\end{Example}

\msk
The weakly complete unital algebras form a category 
$\cal{W\hskip-2pt A}$ in a
natural way with morphisms preserving multiplication and identities in
the usual way.

\msk

If $\cal A$ and $\cal B$ are commutative monoidal categories, 
then a functor $F\colon{\cal A}\to{\cal B}$ is a {\it functor
of monoidal categories} or, more shortly{\color{\red},}
a {\it multiplicative}
functor{\color{\red},} if it firstly induces a functor on the 
underlying categories
and if it secondly respects multiplication and identities so that
there are natural isomorphisms 
$$F(A_1\otimes_{\cal A}A_2)\cong F(A_1)\otimes_{\cal B}F(A_2),\hbox{ and}$$
$$ F(\IdObject_{\cal A})\cong \IdObject_{\cal B}.$$

It is then a routine to prove the following fact:
\msk

\begin{Proposition} \label{p:monoids} \quad {\rm(a)} A multiplicative functor
$F\colon {\cal A}\to{\cal B}$ between commutative monoidal
categories maps monoids in $\cal A$ to monoids in $\cal B$
and induces a functor between the categories of monoids in
the respective categories.

\msk

\nin{\rm(b)}\quad In particular, if $F$ implements an equivalence of
categories {\rm (see \cite{book3}, Definition A3.39)}, then 
the respective categories of monoids 
are likewise equivalent. 
\end{Proposition}

\msk

Let us consider the two categories
$\cal W$ and $\cal V$ and denote by ${\cal V}^{\rm op}$
the opposite category in which all arrows are formally reversed
(cf.\ \cite{book3}, Definition A3.25 and  paragraph preceding it).
Then by Theorem \ref{th:vect}, the functors
$$E\mapsto E^*:{\cal V}^{\rm op}\to {\cal W},\hbox{ respectively, }
V\mapsto V':{\cal W}\to {\cal V}^{\rm op}$$
implement an equivalence of categories, and these functors
are multiplicative by Proposition \ref{p:multi}. 
Hence by Proposition \ref{p:monoids} we have
the following lemma:

\msk 

\begin{Lemma} \label{l:predual} \quad {\it The category of 
weakly complete unital 
topological algebras is equivalent to the category of monoids
in ${\cal V}^{\rm op}$.}
\end{Lemma}

\msk

 But what is a monoid in ${\cal V}^{\rm op}$?
\msk
\begin{Definition} \label{d:comonoids} \quad Let $\cal A$ 
be a commutative monoidal 
category. Then a {\it co\-monoid} in $\cal A$ is an object $C$
together with morphisms $c\colon C\to C\otimes C$ and 
$k\colon C\to \IdObject$ such that 
\def\mdvertspace{\vphantom{A_{A_A}}}\def\muvertspace{\vphantom{A^{A^A}}}

$$\begin{matrix}%
C\ot(C\ot C)&\mapright{\alpha_{CCC}}&(C\ot C)\ot C\\
\lmapup{\id_C\ot c}&&\mapup{c\otimes \id_C}\\
C\otimes C && C\otimes C\\
\lmapup{c}&&\mapup{c}\\
C&\lmapright{\id_c}&C.\\
\end{matrix}
$$

\medskip
 
\cen{and}

\medskip 

$$
\begin{matrix}%
\mdvertspace \IdObject\ot C&\mapleft{k\ot\id_C}&C\ot C&%
         \mapright{\id_C\ot k}&C\ot \IdObject\\
\lmapup{\iota_C}&&\lmapup{c}&&\mapup{\iota'_C}\\
\muvertspace C&\lmapleft{\id_C}&C&\lmapright{\id_C}&\hphantom{.}C.\\
\end{matrix}
$$
\medskip

\nin are commutative (for natural isomorphisms 
$\iota_C\colon C\to \IdObject\otimes C$
and $\iota_C'\colon C\to C\otimes \IdObject$.)
\end{Definition}

An almost trivial example of comonoids are obtained in any category
${\cal C}$ with finite products and a terminal object $\IdObject$ (that is,
an object such that $\Hom_{\cal C}(X,\IdObject)$ is singleton for all objects
$X$: see \cite{book3}, Definition A3.6). Indeed there is a unique 
``diagonal'' morphism $c_X:X\to X\times X$, and a unique morphism
$k_X\colon X\to \IdObject$ and these two morphisms make {\it every} object $X$
endowed with $c_X$ and $k_X$ a coassociative comonoid.
This applies to the category of sets with the cartesian product
and constant functions to the singleton set. (The verification is
a straightforward exercise.)

\msk
From the definition it is clear that a comonoid in $\cal A$ is
exactly a monoid in the opposite category ${\cal A}^{\rm op}$.
If $\cal A$ is a category of vector spaces over some field, notably
in the case of {\color{\red}$\R$ or $\C$}, another name is present in the literature
(see \cite{mich}):
\msk
\begin{Definition} \label{d:coalgebra} \quad A comonoid in 
the commutative monoidal
category $\cal V$ of $\K$-vector spaces is called a 
{\it coalgebra} over $\K$. 
\end{Definition}

\msk

Accordingly, Lemma \ref{l:predual} translates into the following 
statement: 

\msk

\def\W{\cal W}
\def\WA{\cal W\hskip-2pt A}
\def\CA{\cal C\hskip-2pt A}

\begin{Theorem} \label{th:predual} \quad{\it The category   
$\WA$ of weakly complete
unital topological $\K$-algebras is dual to the category of
 $\K$-coalgebras $\CA$.} 
\end{Theorem}

\msk
These observations conclude our collection of preliminary
concepts. Theorem \ref{th:predual} is all but profound.
So let us pay attention where it may lead us.

        \section{Weakly Complete Unital Algebras and
                      their Group of Units}

There is a dominant theorem in the theory of
coalgebras, called the Fundamental Theorem of Coalgebras
(see \cite{mich}, Theorem 4.12, p. 742), attributed to
{\sc Cartier}.
For us, the following version is relevant. It should be clear that
a vector subspace $S$ of a coalgebra $C$ is a subcoalgebra if 
$c_C(S)\subseteq S\otimes S$ and $k_C(S)=\R$. 

\msk

\begin{Theorem} {\rm (Fundamental Theorem of Coalgebras)} \quad 
\label{th:cartier} Every 
coalgebra $C$
is the directed union of the set
of its finite dimensional subcoalgebras. 
\end{Theorem}

\msk

This is sometimes formulated as follows: {\it Every coalgebra
is the injective limit of its finite dimensional subcoalgebras.} 

Now if we take Theorems \ref{th:predual} and \ref{th:cartier}
together, we arrive at the 
following theorem  \cite{dah}. Its 
consequences, as we shall see, are surprising.
(For a discussion of limits in the sense of category theory
see \cite{book3}, Definition A3.41 ff., and in the concrete case
of limits and notably projective limits of
topological groups see \cite{probook}, pp.~63ff., respectively, pp.~77ff.)

A projective limit of of topological groups is {\it strict}
if all bonding morphisms and all limit morphisms are surjective
(see \cite{book3}, 1.32 or \cite{probook}, Definition 1.24).
We shall call a projective limit of topological groups
{\it a strict projective limit of quotients} 
if all bonding maps and all
limit morphisms are surjective {\it and open}; that is,
are quotient morphisms.
In the situation of vector spaces, Theorem A2.12 
of \cite{probook} implies that for 
an injective morphism $f\colon E_1\to E_2$
in the category $\cal V$ of vector spaces, the dual morphism

\cen{$f^*\colon E_2^*\to E_1^*$} 

\nin in the category $\cal W$ of
weakly complete vector spaces is automatically surjective and open. 

\msk

\begin{Theorem} \label{th:fundamental}
{\rm (The Fundamental Theorem of Weakly Complete
Topological Algebras)\quad} Every weakly complete unital 
topological $\K$-algebra is the strict projective limit of
{\color{\red} a projective system of quotient morphisms 
between}   its finite
dimensional unital quotient-algebras.
\end{Theorem} 

\msk 

\nin The literature on locally compact groups shows considerable attention
to structural results derived from the information that
a group, say, is a projective limit of Lie groups. It is therefore
remarkable that a result concluding the presence of a projective
limit of additive Lie groups emerges out of the
vector space duality between $\cal V$ and $\cal W$ and the Fundamental
Theorem on Coalgebras.
 
\msk
For a weakly complete topological unital algebra $A$ let
$\I(A)$ denote the filter basis of closed two-sided ideals $I$ of $A$
such that $\dim A/I<\infty$. We apply Theorem 1.30 of \cite{probook}
and formulate:

\msk

\begin{Corollary} \label{c:ideal-filter} In a weakly complete 
topological unital algebra $A$ each neighborhood
of $0$ contains an ideal $J\in\I(A)$. That is, the filter basis $\I(A)$
converges to 0. In short{\color{\red},} 
$\lim \I(A)=0$ and $A\cong\lim_{J\in\I(A)} A/J$.
\end{Corollary}

\msk

\nin If $A$ is an {\it arbitrary} unital $\K$-algebra and $\I(A)$ is the lattice
of all of its two-sided ideals $J$ such that $\dim A/J<\infty$, then
$\lim_{J\in\I(A)} A/J$ is a weakly complete unital algebra.
It is instructive to apply this procedure 
to the polynomial algebra ${\color{\red} A=}\K[X]$ in one variable 
and obtain the weakly complete algebra
$$\K\<X\>\defi \lim_{J\in\I(\K[X])} \K[X]/J,$$ which
we claim to satisfy the following
universal property: 

\msk

\begin{Lemma} \label{l:universal} 
For any element $a$ in a weakly complete $\K$-algebra $A$ there
is a unique morphism of weakly complete algebras 
$\phi\colon \K\<X\>\to A$ such that $\phi(X)=a$.
\end{Lemma}

\begin{Proof} As a first step, 
let us observe that the algebra generated by 
$X$ in $\K\<X\>=\lim_{J\in\I(\K[X])} \K[X]/J$ is a dense subset. 
This implies that the map $\phi$ in the statement of the Lemma 
is unique if it exists. It remains to show its existence.

Recall that $A\cong\lim_{J\in\I(A)} A/J$. 
Note that for $J\in\I(A)$ there is a unique morphism 
$\phi_J\colon\K[X]\to A/J$ sending a polynomial $p$ to $p(a)+J$.
Let us look at the following diagram:
$$
\begin{matrix}
\K[X]&\mapright{\inc}&\K\<X\>&\mapright{\lambda}&A\\
\lmapdown{\id}&&\mapdown{\rho_J}&&\mapdown{\quot}\\
\K[X]&\lmapright{\quot}&{\frac{\K[X]}{\ker\phi_J}}&\lmapright{\phi_J'}&A/J,\\
\end{matrix}
$$

\bsk
\noindent where $\lambda$ is the map we want to define, 
where $\phi'_J$ is 
induced by $\phi_J$, and where $\rho_J$ is the limit map from
$\K\<X\>=\lim_{I\in\I(\K[X])}\K[X]/I$. Define 
$\lambda_J=\phi'_J\circ \rho_J$. The ``diagonal'' morphisms
$\lambda_J\colon \K\<X\>\to A/J$, $J\in\I(A)$ are seen 
to be compatible
with the morphisms $A/J_2\to A/J_1$ for $J_1\supseteq J_2$ in $\I(A)$.
Accordingly, by the universal property of the limit 
(cf.~\cite{book3}, Definition A3.41), the unique fill-in morphism 
$\lambda$ exists as asserted.
\end{Proof}

\nin\quad Let $\P$ denote the set of all 
irreducible polynomials $p$ with leading coefficient $1$.
Then we have the following:

\msk

\begin{Lemma} \label{l:3.12} 
There is an isomorphism of weakly complete
$\K$-algebras
$$ \K\<X\>\cong \prod_{p\in\P}\K_p\<x\>,\mbox{ where } 
\K_p\<X\>=\lim_{k\in\N}{\frac{\K[X]}{(p^k)}}.$$
\end{Lemma}

\begin{Proof} We recall that every 
ideal $J\in\I(\K[X])$ is generated by 
a nonzero polynomial $f=f(X)$, that is $J=(f)$, 
since $\K[X]$ is a principal ideal domain. 
Furthermore, each polynomial $f$ admits a unique 
decomposition into irreducible factors:

$$
\I(\K[X])=\left\{\left(\prod_{p\in\P} p^{k_p}\right) : 
     (k_p)_{p\in\P}\in(\N_0)^{(\P)}\right\}.
$$
Here, $(\N_0)^{(\P)}$ denotes the set of all families of nonnegative 
integers where all but finitely many indices are zero.
For each $f=\prod_{p\in\P} p^{k_p}$ we have 
$$
\K[X]/(f)\cong\prod_{p\in\P}\K[X]/(p^{k_p})
$$ 
by the Chinese Remainder Theorem.

This enables us to rewrite the projective limit in the definition 
of $\K\<x\>$ as

$$
 \lim_{J\in\I( \K[X])} \K[X]/J\to%
\prod_{p\in\P}\left(\lim_{k\in\N}{\frac{\K[X]}{(p^k)}}\right).
$$
\vskip-32.5pt
\end{Proof}

\msk

\nin
We remark that if $p\in\P$ is of degree $1$, 
the algebra $\K_p\<X\>$ is isomorphic to $\K[[X]]$, 
the power series algebra in one variable.

Since for $\K=\C$, all $p\in\P$ are of degree $1$, 
it follows that the algebra $\C\<X\>$ is isomorphic 
to $\C[[X]]^\C$.

\msk

In the case $\K=\R$, the polynomials $p\in\P$ are of degree one or two. 
For $p=X-r$ with a number $r\in\R$ the algebra $\R_p\<X\>$ is isomorphic to $\R[[X]]$. 
But for polynomials $p\in\P$ of degree two, 
the situation becomes more complicated.

\bsk

\subsection{The Group of Units: Density} 
An element $a$ in an algebra $A$ is called a {\it unit}
 if it has a multiplicative inverse, that is, there exists an
element $a'\in A$ such that $aa'=a'a=1$. The set $A^{-1}$
 of units of an algebra is a group with respect to 
multiplication.

\msk

\begin{Lemma}  \label{l:inversion-continuous}The group of units 
$A^{-1}$ of a weakly complete
unital algebra $A$ is a topological group.
\end{Lemma}

\begin{Proof} We must show that the function 
$a\mapsto a^{-1}: A^{-1}\to A^{-1}$ is continuous.
In every finite dimensional real or complex unital algebra,
the group of units is a topologial group. 
This applies to each factor algebra $A/I$, for $I\in\I(A)$.
Then $a\mapsto a^{-1}I: A^{-1}\to (A/I)^{-1}$ is 
continuous for all $I\in \I(A)$. Since the isomorphism 
$A\cong \lim_{I\in\I(A)} A/I$ holds
also in the category of topological spaces, the continuity of
$a\mapsto a^{-1}$ follows by the universal property of the 
limit (see \cite{book3}, Definition A3.41).
\end{Proof}
\msk

We remark that there exist topological algebras in which inversion in
the group of units is discontinuous. 
(See e.g.\ \cite{dahii}, Example 3.12.)

\bsk

The prescription which assigns to a weakly complete unital algebra
$A$ its group of units $A^{-1}$ is a functor from the category of 
weakly complete unital algebras to the category of topological groups.
This functor preserves products and intersections, hence arbitrary limits
(see e.g.\ \cite{book3}, Proposition A3.53). Thus $A\cong\lim_{J\in\I(A)} A/J$
implies $A^{-1}\cong\lim_{J\in\I(A)} (A/J)^{-1}$.
Since the group of units of a finite dimensional 
unital algebra is a (finite dimensional) linear Lie group 
(see \cite{book3}, Definition 5.32) we have
\msk
\begin{Lemma} \label{l:4.2} The group of units 
$A^{-1}$ of a 
weakly complete unital real or complex algebra $A$ is a projective limit of 
linear Lie groups.
\end{Lemma}
\msk

\nin Due to a Theorem of Mostert and Shields \cite{most}, the group of units of
a topological monoid on a locally euclidean space has an open 
group of units. 
 This applies, in particular, to the multiplicative semigroup $(A,\cdot)$
of any 
finite dimensional real or complex algebra $A$. Howeveer, in this case, one has more elementary
linear algebra arguments to be aware of this fact. Indeed, let $(A,+)$ 
denote the vector space underlying $A$ and
$a\mapsto L_a:A\to \Hom((A,+),(A,+))$ the representation of $A$ into the
algebra of all vector space endomorphisms of $(A,+)$ given by
$L_a(x)=ax$. If we set $\delta(a)=\det(L_a)$, then we obtain a 
morphism of multiplicative monoids 
$\delta\colon (A,\cdot)\to (\R,\cdot)$
in such a fashion that $A\setminus A^{-1}=\delta^{-1}(\{0\})$. This set is a closed
nowhere dense algebraic hypersurface in $A$. Thus we have
\msk

\begin{Lemma} \label{l:4.3} 
If $A$ is a finite dimensional real or complex unital algebra,
then the group $A^{-1}$ of units is a dense open subgroup of 
the monoid $(A,\cdot)$. 
\end{Lemma}

\nin It may be helpful to consider some examples:

\msk

\begin{Example} (a) In the 
weakly complete algebra $A=\R^\N$ of
Example \ref{e:3.2}.1a,
the identity $1=(e_n)_{n\in\N}$, $e_n=1\in\R$ is the limit of the
sequence of $(a_m)_{m\in\N}$, $a_m=(a_{mn})_{n\in\N}$ of nonunits, where
$$a_{mn}=\begin{cases}0, &\hbox{if $n=m$},\\ 1, &\hbox{otherwise}.\\
           \end{cases}
$$
Hence $A^{-1}$ fails to be open in $A$ while it is dense in $A$.

\msk

(b) In the Examples \ref{e:3.2}.2 and \ref{e:3.2}.3
 of weakly complete unital algebras  
$A$, we have a maximal ideal $M$ (in case (2) containing all elements
with $a_0=0$ and in case (3) equalling $V\times\{0\}$) such that
 $A^{-1}=A\setminus M$. Thus in these cases $A^{-1}$ is open 
(regardless of dimension).
\end{Example}

In order to establish
a first important result on the group of units $A^{-1}$ of a weakly
complete unital algebra $A$, namely, its density in $A$, we need to prove 
it in the special case $$A=\K_p\<X\>=\lim_{k\in \N}\K[X]/(p^k)\leqno(\#)$$ 
of Lemma 3.12.
\msk

\begin{Lemma} \label{l:density} {\rm (Density Lemma)}\quad 
For each irreducible polynomial $p$ over 
$\K$ with leading coefficient $1$, the weakly complete algebra 
$A\defi\K_p\<X\>$ is a local ring and its
group $A^{-1}$ of units 
is open and dense in $A$. 
\end{Lemma}

\begin{Proof} 
Let $\pi \colon A\to \K[X]/(p)$ denote the limit morphism for $k=1$
in $(\#)$ and let $J\defi\ker \pi$. 
For every $f\in J$, the series $\sum_{m=0}^\infty f^m$ converges in 
$A$ to $(1-f)^{-1}$. So 
$$1-J\subseteq A^{-1}.\leqno(1)$$ 

\nin
Now let $f\in A\setminus J$. Since $F\defi \K[X]/(p)$ is a field, 
$\pi(f)$ has an inverse in $F$. Thus there is an element $g\in A$
 with $h\defi fg\in 1-J$. By (1) $h^{-1}$ exists and
$fgh^{-1}=1$. Hence $f$ is invertible. This shows that
$$A\setminus J\subseteq A^{-1}.\leqno(2)$$
Trivially $A^{-1}\cap J=\emptyset$ and so equality holds in (2).

This shows that the closed ideal $J$ is maximal and thus
$A$ is a local ring. Moreover, 
$A^{-1}=A\setminus J=\pi^{-1}(F\setminus\{0\})$
is open and dense as the inverse of a dense set under an open
surjective map.
\end{Proof}

After the Density Lemma \ref{l:density} we notice that 
for an irreducible polynomial $p$ the algebra $\R_p\<X\>$
is a local ring with maximal ideal $p\R_p\<X\>$ such that

$$\frac{\R_p\<X\>}{p\R_p\<X\>}\cong
\begin{cases}\R &\mbox{if }\deg p=1,\\
             \C &\mbox{if }\deg p=2.\\
 \end{cases}$$
For $p=p(X)=X$ we have $\R_p\<X\>=\R[[X]]$ as in Example 2.2.2.
For $J=\R$, Example 2.2.1 can be obtained as a quotient of $\R\<X\>$, 
firstly, by taking for each $r\in\R$ the irreducible polynomial $p=X-r$,
secondly, by noticing that $\R\cong\frac{\R[X]}{X-r}=\frac{\R[X]}p$ is
a quotient of $\frac{\R[X]}{p^k}$ for each $k\in\N$ and thus
of $\lim_{k\in\N}\frac{\R[X]}{p^k}$, and, finally, by 
 determining a quotient morphism
$$\R\<X\>\cong\prod_{p\in\P} \lim_{k\in\N}\frac{\R[X]}{(p^k)}
\to \prod_{r\in\R} \frac{\R[X]}{(X-r)}\cong \R^\R.$$

\nin In passing, we note here again the considerable ``size'' 
of $\K\<X\>$.
\bsk

\begin{Theorem} \label{th:density}{\rm (The First Fundamental Theorem 
on the Group of Units)}
For any weakly complete unital $\K$-algebra $A$,
the group $A^{-1}$ of units is dense in $A$.
\end{Theorem}

\begin{Proof} Let $0\ne a\in A$ and let $V$ denote an open neighborhood of 
$a$ of $A$. According to Lemma \ref{l:universal} there is a morphism
 $\phi\colon \K\<X\>\to A$ with $\phi(X)=a$. Then $U\defi \phi^{-1}(V)$ 
is an open neighborhood of $X$ in $\K\<X\>$. If we find a unit
$u\in \K\<X\>^{-1}$ in $U$, then $\phi(u)\in V\cap A^{-1}$ is a unit,
and this will prove the density of $A^{-1}$ in $A$.
By Lemma \ref{l:3.12} we have $\K\<X\>\cong \prod_{p\in\P}\K_p\<X\>$,
and so the problem reduces to finding a unit near $X$ in 
$\K_p\<X\>$ 
for each $p\in\P$. The preceding Density Lemma \ref{l:density}
says that this is possible.
\end{Proof}

\bsk
\subsection{The Exponential Function}
Every finite dimensional unital $\K$-algebra is, in particular, a
unital Banach algebra over $\K$ with respect to a suitable norm. 
By \cite{book3}, Proposition 1.4, in any unital Banach algebra $A$ 
over $\K$ the group
$A^{-1}$ of units is an open subgroup of the monoid $(A,\cdot)$, and 
 it is a (real) linear Lie group with 
Lie algebra $\L(A)=A_{\rm Lie}$, the real 
vector space underlying $A$ with the 
Lie bracket given by $[x,y]=xy-yx$, with the exponential function
$\exp\colon \L(A^{-1})\to A^{-1}$ given by the
everywhere absolutely convergent power series
$\exp x=\sum_{n=0}^\infty\frac1{n!}\.x^n$. 
(For $\K=\R$ this is discussed extensively
in \cite{book3}, Chapter 5, notably
\cite{book3}, Definition 5.32.)

\msk

Now let $A$ be a weakly complete unital $\K$-algebra. Every closed
(2-sided) ideal $J$ of $A$ is a closed Lie algebra ideal
of $A_{\rm Lie}$. We apply the Theorem \ref{th:fundamental}
and note that
the Lie algebra $A_{\rm Lie}$ is (up to natural isomorphism
of topological Lie algebras) the strict projective limit of quotients
$$\lim_{J\in\I(A)} \left(\frac{A}{J}\right)_{\rm Lie}\subseteq 
\prod_{J\in\I(A)}\left(\frac{A}{J}\right)_{\rm Lie}$$ 
of its finite dimensional quotient algebras 
and therefore is a pro-Lie algebra. Each of 
these quotient Lie algebras is the domain of an exponential function 
$$\exp_{A/J}\colon A_{\rm Lie}/J\to (A/J)^{-1}\subseteq A/J,
\quad (\forall a_J\in A/J)\,
 \exp_{A/J} a_J=\sum_{n=0}^\infty\frac1{n!}\.a_J^n.$$
This yields a componentwise exponential function on
$\prod_{J\in\I(A)}A/J$ which respects the bonding morphisms of 
the subalgebra $\lim_{J\in\I(A)}A/J$. Thus we obtain the following
basic result which one finds in \cite{dah, hofmor}.
\msk

\begin{Theorem} {\rm(The Second Fundamental Theorem on the Group of Units)}
\label{th:second}
   If $A$ is a weakly complete unital $\K$-algebra,
then the exponential series $1+a+\frac1{2!}a^2+\cdots$ 
converges on all of $A$ and defines the exponential function
$$\exp_A \colon A_{\rm Lie}\to A^{-1},\quad 
        \exp_A a=\sum_{n=0}^\infty \frac1{n!}a^n$$
of the pro-Lie group $A^{-1}$. The Lie algebra $\L(A^{-1})$
of the pro-Lie group $A^{-1}$ may be identified with 
the topological Lie algebra
 $A_{\rm Lie}$, whose underlying weakly
complete vector space is the underlying weakly complete 
vector space of $A$.
\end{Theorem}

\msk
It is instructive to observe that while we saw any weakly 
complete associative unital algebra to be a projective
limit of finite dimensional quotient algebras, it would be
incorrect to suspect that every weakly complete real Lie algebra
was a projective limit of finite dimensional quotient algebras.
Indeed consider the weakly complete vector space $L\defi\R^\N\times R$
with the Lie bracket 
$[\((x_n)_{n\in\N},s\), \((y_n)_{n\in\N},t\)]
=\(s\.(y_{n+1})_{n\\in\N} -t\.(x_{n+1})_{n\in\N},0)\)$. 
The $L$ is a weakly complete Lie algebra which can be
shown to have no arbitrarily small cofinite dimensional ideals
and so cannot be the projective limit of finite dimensional
quotients.

\subsection{The Group of Units of Weakly Complete Algebras}
For the investigation of the group of units of a weakly complete
unital algebra, this opens up the entire 
Lie theory of pro-Lie groups $G$, for which in the case of $\K=\R$
 we refer to \cite{probook}.
For instance, Lemma 3.29 of \cite{probook} tells us that $A^{-1}$
is a pro-Lie group.

\msk

\noindent{\bf The units in finite dimensional real or complex algebras.}
For a finite dimensional unital algebra $A$ let $G\defi A^{-1}$ the
subset of units, i.e. multiplicatively invertible elements. the set of
pairs $(a,b)\in A\times A$ such that $ab=1$ where $1$ is the identity
of $A$ is the set of solutions of a finite sequence of 
equations over $\K$ and is the graph of the function 
$$a\mapsto a^{-1}:G \to G$$
and therefore is a  algebraic variety over $\K$. It is homeomorphic to
$G$ in the topology induced from the unique $\K$-vector space
topology of $A$. For $\K=\R$ and $\K=\C$ this identifies $G$
as a real Lie group. This raises the question about information
on $G/G_0$ where $G_0$ is the identity component of $G$.
 {\color{\red}Recall that a topological group $G$ is called
{\it almost connected} if the factor group $G/G_0$ is compact.}

\begin{Lemma} \label{l:4.6} 
For a finite dimensional real unital algebra $A$
the group of units $A^{-1}$ has finitely many components and thus is
almost connected. 
\end{Lemma}

\begin{Proof} See \cite{boch}, Theorem 2.1.5.
\end{Proof}
\msk
\nin If $\K=\R$ is the field of real numbers and $A=\R^n$, then 
$A^{-1}/(A^{-1})_0\cong \Z(2)^n$ has $2^n$ components.

\bsk

\begin{Lemma} \label{l:canslemma} {\rm {\color{\red}(}Mahir 
Can's Lemma{\color{\red})}} 
For a finite dimensional complex unital algebra $A$
the group of units $A^{-1}$ is connected.
\end{Lemma}

\begin{Proof}  We write $G=A^{-1}$ for the group of units of $A$
and let $N$ denote the nilradical of the algebra $A$. Then for some
natural number $k\ge2$ we have $N^k=\{0\}$ and so for each $x\in N$
we have $(1-x)^{-1}=1+x+\cdots+x^{k-1}$. Thus $1+N$ is a subgroup
of $G$ which is immediately seen to be connected. 
Since $N$ is an ideal, for any $g\in G$ we have
$N=gg^{-1}N\subseteq gN\subseteq N$ and so $gN=N$ and
$g(1+N)=g+N\in A/N$, that is,  $G/(1+N)\subseteq A/N$ for
the semisimple  factor algebra $A/N$ of $A$ modulo its
nilradical.

\msk

By a classical theorem of {\color{\red}Wedderburn} there is an $m$-tuple
$(n_1,\dots,n_m)$ of nonnegative integers such that
$$A/N\cong M_{n_1}(\C)\oplus\cdots\oplus M_{n_m}(\C),$$
where $M_n(\C)$ denotes the $n\times n$ complex
matrix algebra ($\cong \Hom(\C^n,\C^n)$).
(See e.g.~Corollary 2.66 in \cite{bresar}).
The group ${\rm GL}(n,\C)=M_n(\C)^{-1}$ of units in the full matrix 
algebra is connected, and so the group of units
$H\defi (A/N)^{-1}$ of the algebra $A/N$  is connected.
(See e.g.~\cite{good}, Theorem 2.2.5.)

\msk

Let $p\colon A\to A/N$ be the quotient map. We saw
that $p$ induces the morphism $p|G\colon G\to H$ 
{\color{\red} of}
complex Lie groups. By Lemma \ref{l:4.3}, $G$ is
dense and open in $A$ and $H$ is open in $A/N$.
The function $p$ is open, and thus the function
$p|G$ is an open morphism of topological groups
so that $p(G)$ is open, hence closed in $H$.
The density of $G$ in $A$ implies the density
of $p(G)$ in $p(A)=A/N$ and so in $H$. Therefore
$p(G)=H$. Thus we have an exact sequence
$$1\to 1{+}N \to G\mapright{p|G}H\to1.$$
So, since both $1+N$ and $H$ are connected,
$G$ is connected which we had to show.
\end{Proof}

\msk

\noindent{\bf Quotient morphisms of algebras.}
Let $f\colon A\to B$ be a surjective morphism
of weakly complete unital algebras.
The restriction and corestriction of this function
$f|A^{-1}:A^{-1}\to B^{-1}$ is
a morphism of pro-Lie groups. 

Then $\L(f|A^{-1})\colon\L(A^{-1})\to\L(B^{-1})$ 
is $f$ itself considered as a surjective morphism of 
pro-Lie algebras.
Then by \cite{probook}, Corollary 4.22.iii, we have
$$f(\<\exp_A(A)\>)=(f|A^{-1})\<\exp_A(A)\>=\<\exp_B(B)\>. \leqno(1)$$

Assume now that $\dim B<\infty$.
Then $B^{-1}$ is a Lie group and so
$$\<\exp_B(B)\>=(B^{-1})_0.\leqno(2)$$ 
Thus $ (B^{-1})_0\subseteq f(A^{-1})$
and so 
$$f(A^{-1})\hbox{ is open and closed in }B^{-1}.\leqno(3)$$

By Theorem \ref{th:density}, the group of units
$A^{-1}$ is dense in $A$, whence
$f(A^{-1})$ is dense in $f(A)=B$, and so, in particular, 
$f(A^{-1})\subseteq B^{-1}$ is dense in $B^{-1}$. 
But $f(A^{-1})$ is also closed in $B^{-1}$ by (3),
and thus $f(A^{-1})=B^{-1}$.

\msk

Let $I=\ker f$. Then $A^{-1}\cap(1+I)=\ker f|A^{-1}$.
Thus we have a unique bijective morphism 
$$f'\colon \frac{A^{-1}}{A^{-1}\cap(1+I)}\to B^{-1}$$
such that 
$$(\forall x\in A^{-1})\, f'(x(A^{-1}\cap(1+I)))=f(x).$$

Since $A/I\cong B$, as vector spaces, we know that 
$A^{-1}/(A^{-1}\cap(1+I)) \cong(A/I)^{-1}$ is a Lie 
group with finitely many components by Lemma \ref{l:4.6}, 
which is, therefore,
$\sigma$-compact. Hence $f'$ is
an isomorphism by the Open Mapping Theorem for 
Locally Compact Groups (see e.g.\ \cite{book3},
Exercise EA1.21 on p.704). 
If $q\colon A^{-1}\to A^{-1}/\ker(f|A^{-1})$ is the quotient
morphism, then $f|A^{-1}=f'\circ q$. Hence $f|A^{-1}$ is
an open morphism, and we have shown the following 

\begin{Lemma} \label{l:4.8} Let
$f\colon A\to B$ be a surjective morphism
of weakly complete unital algebras and assume that $\dim B<\infty$. 
Then $f$ induces
a quotient morphism $f|A^{-1}:A^{-1}\to B^{-1}$ 
from the pro-Lie group $A^{-1}$ onto the Lie group $B^{-1}$.
\end{Lemma}

Keep in mind that a quotient morphism is an open map!
We also note that Lemma \ref{l:4.8} remains true for $\K=\C$.
\bsk

\nin We apply this to any weakly complete unital algebra $A$ and
any finite dimensional quotient algebras $B=A/J$ for $J\in\I(A)$
and thus sharpen Lemma \ref{l:4.2} in a significant way:

\begin{Theorem} \label{th:strict-quotients}
{\rm(The Third Fundamental Theorem on Units)}\quad For 
a weakly complete unital $\K$-algebra $A$
let $\I(A)$ be the set of two-sided closed ideals $J$ such that 
$\dim A/J<\infty$. Then $\I(A)$ is a filter basis converging to $0$
and the group of units $A^{-1}$ is a {\rm strict} projective limit 
of quotients of linear almost connected
Lie groups $A^{-1}/(A^{-1}\cap(1+J))$ isomorphic to 
$(A/J)^{-1}$ via the map $g(A^{-1}\cap(1+J))\mapsto g+J$
{\color{\red} as $J$ ranges through $\I(A)$,}  
such that each limit morphism agrees with the natural 
quotient morphism.
\end{Theorem}

\nin 
With this theorem, at the latest, it becomes clear 
that we had to resort to the
structure theory of pro-Lie groups with almost connected
Lie groups quotients. 
In \cite{almost} the extensive structure and Lie theory of
{\it connected} pro-Lie groups of \cite{probook} was 
lifted to {\it almost connected} pro-Lie group. Yet
it was still unknown whether a projective limit of
almost connected finite dimensional Lie groups was
complete and is, therefore, an almost connected pro-Lie group.
However, this is now clear with the result recorded in
\cite{short}.

\section{Applying pro-Lie Group Theory to the Group of Units
of a Weakly Complete Algebra}

At this point we continue the theory of weakly complete algebras
and prove the following result:

\begin{Theorem} {\rm (The Fourth Fundamental Theorem on the Group 
of Units)} \label{th:fourth}
 In any weakly complete  unital algebra $A$,
the multiplicative group $A^{-1}$ of invertible elements
is almost connected if the ground field is real and
is connected if the ground field is complex. 
\end{Theorem}

\begin{Proof} (a) Let us abbreviate $G\defi {\color{\red}A}^{-1}$.
 Assume that the ground field $\K$ is $\R$.
By Theorem \ref{th:strict-quotients},
 $A=\lim_{J\in\I(A)} A/J$ where $\dim A/J<\infty$
 and $G=\lim_{J\in\I(A)} (A/J)^{-1}$
with finite dimensional Lie groups $(A/J)^{-1}$. By Lemma \ref{l:4.6}, 
the group $(A/J)^{-1}$ is almost connected.
Let $L=G/N$ be any Lie group quotient of
$G$. Let $U$ be an identity neighborhood in $G$ 
such that $UN=U$ and
that $U/N$ has no nonsingleton subgroup. Since $\lim\I(A)=0$
we have $\lim_{J\in \I(A)} 1+J=1$. So there is a 
$J\in\I(A)$ such that $1+J\subseteq U$. Since 
 $G\cap(1+J)$ is a multiplicative
subgroup, $G\cap(1+J)\subseteq N$. Hence $L$ is a homomorphic
image of $G/(G\cap(1+J))\cong(A/J)^{-1}$. This group has finitely
many components,
and so $L$ is almost connected.
Now Theorem 1.1 of \cite{short}
 shows that $G$ is almost connected in the real case. This means
that $G/G_0$ is a compact totally disconnected group and thus has
arbitrarily small open-closed subgroups.

(b) Assume now $\K=\C$. In particular, $A$ is a real weakly complete
unital algebra and thus $G/G_0$ has arbitrarily small open-closed subgroups.
Suppose that $G\ne G_0$. Then there is a proper open subgroup $H\subseteq G$.

According to Theorem \ref{th:strict-quotients}
$G$ is the strict projective limit of the complex
Lie group quotients $G/(G\cap(1+J))=(A/J)^{-1}$,  $J\in\I(A)$.
Thus we find a $J\in\I(A)$ so small  that $G\cap(1+J)\subseteq H$.
Then $H/(G\cap(1+J))$ is a proper open subgroup of $(A/J)^{-1}$.
However this complex Lie group is  connected by Lemma \ref{l:canslemma}. 
This contradiction shows that the assumption $G\ne G_0$ is
false, and thus that $G$ is connected. This completes the 
proof of the theorem.
\end{Proof}

\msk

After this information is secured we invoke basic results 
of pro-Lie groups combined from 
\cite{probook}, 12.81 on p.551, and \cite{almost}.
We let $A$ denote any weakly complete unital algebra
(see Definition \ref{d:basic-definition}) and
denote by $G$ its group of units $A^{-1}$.

\begin{Theorem} \label{th:vector-space-splitting}
The group of units $G$ of a weakly complete unital algebra $A$
 contains a maximal compact subgroup $C$
and $A$ contains up to four closed vector subspaces
$V_1,\dots,V_m$ such that 
$$(c,X_1,\dots,X_m)\mapsto c\exp X_1\cdots\exp X_m\colon
C\times V_1\times\cdots\times V_m\to G$$
is a homeomorphism. Every compact subgroup of $G$
has a conjugate in $C$. The group $G$ is homeomorphic
to a product space $\R^J\times C$ for a suitable set $J$.
\end{Theorem} 

One has more detailed structural information on $G$ when
needed: If $N(G)$ denotes the nilradical of $G_0$
(see \cite{probook}, Definition 10.40 on p.~447),
then for each $1\le k\le m$ the product $N(G)\exp V_k$
is a prosolvable closed subgroup (see Definition 10.12,
p. 424 of \cite{probook}).

\msk
Since $G$ is homeomorphic to a space of the form $\R^J\times C$
for a weakly complete vector space $\R^J$ and a compact group
$C$ we know that all of the algebraic-topological properties
of $G$ are those of a compact group, since $\R^J\times C$ and
$C$ are homotopy equivalent. (Cf.\ \cite{hofmor}.)

\msk

From Lemma 1.1(iii) in \cite{short}, further Theorem \ref{th:vector-space-splitting},
and \cite{book3}, Theorem 9.41 on p.~485, and \cite{book3}, 
Corollary 10.38 on p.~572{\color{\red},} 
we also derive the following facts.

\begin{Corollary} Let $G$ be the group of units of a weakly complete unital
algebra. Then

\begin{enumerate}[\rm(i)]

\item $G$ contains a profinite subgroup $D$ such that $G=G_0D$ while 
     $G_0\cap D$ is normal in G and central in $G_0$, and
\item $G$ contains a compact subspace $\Delta$ such that
      $(c,d)\mapsto cd:G_0\times \Delta\to G$ is a homeomorphism. 
\end{enumerate}
\end{Corollary} 

\bsk

\subsection{Limitations of these Results} While we know from \cite{book3},
Corollary 2.29 on p.~43 that every compact group is contained (up to isomorphy)
in the group of units of a weakly compact unital algebra, there are even 
Lie group{\color{\red}s} of dimension as small as 3 
which cannot have any isomorphic
copies contained in the group of units of a weakly complete algebra.

\msk

\begin{Example} \label{e:nonlinear} Among  
the connected noncompact Lie groups $G$
which are not linear, the following 3--dimensional 
examples are prominent in the literature:

(a) (Garret Birkhoff) Let $N$ be the ``Heisenberg group'' of matrices 
$$[x,y,z]\defi\begin{pmatrix}1&x&z\\0&1&y\\0&0&1\end{pmatrix},\ x,y,z\in\R$$
and $Z=\{[0,0,n]:n\in\Z\}\cong\Z$ a central cyclic subgroup of $N$.
Then $G=N/Z$ is a 3-dimensional class 2 nilpotent group
 which is not linear as first observed by G. Birkhoff. (See e.g.
\cite{book3}, Example 5.67ff.) The group $G$ is homeomorphic to 
$\R^2\times \Ss^1$.

\ssk

(b) Let $G$ be the universal covering group of the special
linear group ${\rm Sl}(2,\R)$,
the group of 2 by 2 real matrices of determinant 1. Since 
${\rm Sl}(2,\R)$ is homeomorphic to $\R^2\times\Ss^1$, 
the three dimensional Lie group $G$ is homeomorphic to $\R^3$.
The 3-dimensional Lie group $G$ is not linear 
(See e.g.~\cite{hilnee}, Example 9.5.18.)
For a practical parametrization of $G$ see e.g.~\cite{hilg}, 
Theorem V.4.37 on p.~425 and the surrounding discussion.
\end{Example} 

\msk

\nin
We observe that the universal covering group $G$ of
Sl$(2,\R)$ cannot be a subgroup of the group $A^{-1}$ of
units of a weakly complete real unital algebra $A$.
Indeed {\sc Hilgert and Neeb} show in \cite{hilnee}, Example 9.5.18.,
that every linear representation $\pi\colon G\to{\rm Gl}(n,\R)$
factors through $p\colon G\to{\rm Sl}(2,\R)$ with a 
representation $f\colon {\rm Sl}(2,\R)\to {\rm Gl}(n,\R)$
as $\pi=f\circ p$. So if $z$ is one of the two generators 
of the center of $G$ which is isomorphic to $\Z$, then $p(z^2)=1$
and thus $\<z^2\>$ is contained in the kernel of any finite dimensional
linear representation of $G$. From a subgroup $G$ of $A^{-1}$
we would obtain an injective morphism $\gamma\colon G\to A^{-1}$.
Every $J\in\I(A)$ yields a linear representation 
$q_J\colon A^{-1}\to(A/J)^{-1}$ by Theorem 3.6. 
Thus $q_J\circ \gamma:G\to (A/J)^{-1}$ is a linear representation
which will annihilate $\<z^2\>$ for all $J\in\I(A)$. Thus the
injective morphism $\gamma$ would annihilate $z^2$ which is
a contradiction. 

We leave it as an exercise to show that the group $G$ of 
Birkhoff's Example
4.15(a) cannot be a subgroup of any $A^{-1}$ of a weakly complete
unital algebra 1. (Cf.~\cite{book3}, Example 5.67.)

\bsk

\nin
Still, in the next section we shall show that for each topological group
$G$ there is a ``best possible'' weakly complete unital algebra
$\R[G]$ with a natural ``best possible'' morphism $G\to\R[G]^{-1}$, one
that would be an embedding for every compact group but would fail
to be injective for the universal covering group of Sl($2,\R)$.

\msk

\section{The Weakly Complete Group Algebra 
of a Topological Group}

\msk

Let us complement our discussion at this stage by describing 
a relevant pair of adjoint functors.

So let $\cal W\hskip-2pt A$ be the category of weakly complete unital $\K$-algebras
and $\mathcal{T\hskip-1pt G}$ the category of topological groups. Then 
$$\Omega=(A\mapsto A^{-1}):{\cal W\hskip-2pt A}\to{\cal TG}$$ is a well defined
functor after Lemma \ref{l:inversion-continuous}. 
It is rather directly seen to preserve
arbitrary products and intersections. Hence by \cite{book3}, Proposition A3.51
it is a continuous functor (see loc.cit.\ Definition A3.50), that is, it preserves
arbitrary limits. 

\msk

\subsection{The Solution Set condition} 
In order to conclude from this information that $\Omega$ has in fact a 
left adjoint, we need to verify the so called {\it Solution Set
Condition} (see \cite{book3}, Definition A3.58 on p.~786).

\nin For this purpose we claim that for any topological group 
$G$ in there is a {\it set}
$S(G)$ of pairs $(\phi, A)$ with a continuous group 
 morphisms $\phi\colon G\to A^{-1}$
for some object $A$ of $\cal W\hskip-2pt A$ such that for every pair $(f, B)$,
$f\colon G\to B^{-1}$ with a weakly complete unital algebra $B$ there is
a pair $(\phi, A)$ in $S(G)$ and a $\cal{W\hskip-2pt A}$-embedding
$e\colon A\to B$ such that 
$$f=G\mapright{\phi} A^{-1}\, \mapright{e|A^{-1}} B^{-1},$$ 
where $e|A^{-1}$ denotes the bijective restriction and corestriction of $e$. 

Indeed if $f\colon G\to B^{-1}$ determines a unique smallest 
algebraic unital abstract subalgebra $C$ of $B$ generated by $f(G)$,
then  there is only a set of these ``up to equivalence''.
Then 
on each of these there is only a set of algebra
topologies and, a fortiori, only a set of them for which the 
corestriction is continuous; for each of these, there is at most a set of
algebra completions up to isomorphism.
 So, up to equivalence there
is only a set of pairs $(\phi,A)$, $\phi\colon G\to A^{-1}$ such that the
unital algebra generated by $\phi(G)$ is dense in $A$. Any such set
$S(G)$ will satisfy the claim. 

\msk

Now we are in a position to apply the Left Adjoint Functor Existence
Theorem (see \cite{book3}, Theorem A3.60) to conclude that 
$\Omega\colon {\cal{W\hskip-2pt A}}\to {\cal G}$ has a left adjoint 
$\Lambda\colon {\cal G}\to{\cal{W\hskip-2pt A}}$. We write the weakly complete
unital algebra $\Lambda(G)$ as 
$\K[G]$ and call it the {\it weakly complete group algebra of} $G$.
We summarize this result:

\msk
\begin{Theorem} \label{th:wcga-thm} {\rm(The Weakly Complete Group
Algebra Theorem)} \quad To each topological group $G$ 
there is attached functorially a weakly complete
group algebra $\K[G]$ with a natural morphism $\eta_G\colon G\to \K[G]^{-1}$
such that the following universal property holds:

\noindent
For each weakly complete unital algebra $A$ and each morphism 
of topological groups 
$f\colon G\to A^{-1}$ there exists a unique morphism of 
weakly complete unital algebras \break
$f'\colon \K[G]\to A$ restricting to a morphism 
$f''\colon\K[G]^{-1}\to A^{-1}$ of topological groups
 such that $f= f'\circ\eta_G$.
\end{Theorem}

$$ 
\begin{matrix}&{\rm top\ groups}&&\hbox to 7mm{} &{\rm wc\ algebras}\\
\noalign{\vskip3pt}
\noalign{\hrule}\\
\noalign{\vskip3pt}%
        \mdvertspace G&\mapright{\eta_G}&\K[G]^{-1}&\hbox to 7mm{} &\K[G]\\
\lmapdown{\forall f}&&\mapdown{f''}&\hbox to 7mm{}&\mapdown{\exists! f'}\\
\muvertspace A^{-1}&\lmapright{\id}&A^{-1}&\hbox to 7mm{}&A\\
\end{matrix}
$$
\msk

\nin
{\bf Fact.} {\it 
If $G$ is one of the two groups of {\rm Example \ref{e:nonlinear}}, then the natural
morphism $\eta_G\colon G\to\K[G]$ is not injective.}

However, the adjunction of the functors $A\mapsto A^{-1}$ (on the right)
and $G\mapsto \K[G]$ (on the left) also has a 
back adjunction $$\epsilon_A\colon \K[A^{-1}]\to A$$
such that for each topological group $G$ and each
continuous algebra morphism $f\colon \K[G]\to A$ there is a unique
morphism of topological groups $f'\colon G\to A^{-1}$ such that
$f=\epsilon_A\circ \K[f']$. (Cf.\ \cite{book3}, Proposition A3.36, p.~777).
 The general theory of adjunctions (as e.g. in \cite{book3},
Proposition A3.38, p.~777) now tells us that we may formulate

\msk 
\begin{Corollary} \label{c:alternate} For any weakly complete unital
algebras $A$ and any topological groups $G$ we have

$$ (\forall A)\, \Big(A^{-1}\quad\mapright{\eta_{A^{-1}}}\quad \K[A^{-1}]^{-1}
 \quad\mapright{(\epsilon_A)^{-1}} A^{-1}\Big)
    =\id_{A^{-1}},\quad\hbox{and}$$

\medskip

$$(\forall G)\,\Big(\K[G]\quad\mapright{\K[\eta_G]}\quad \K\big[\K[G]^{-1}\big]
\quad \mapright{\epsilon_{\K[G]}}\quad \K[G]\Big)=\id_{\K[G]}.$$
\end{Corollary}

 In other words: Many topological groups are semidirect factors of unit groups
of weakly complete algebras, for instance if they are unit groups to begin with,
and many weakly complete unital algebras are homomorphic retracts
(semidirect summands) of weakly complete group algebras, for instance,
if they are group algebras to begin with. 

\bsk

However, more importantly, the universal property implies
conclusions which are relevant for the concrete structure theory
of $\K[G]$:

\msk

\begin{Proposition} \label{p:generation} 
	Let $G$ be a topological group.
	Then the subalgebra linearly spanned by $\eta_G(G)$
	in $\K[G]$ is dense in $\K[G]$.
\end{Proposition}

\begin{Proof} 
\ Let $S=\overline{\Span}\(\eta_G(G)\)\subseteq \K[G]$ 
	be the closed subalgebra linearly spanned by $\eta_G(G)$.
	Let $f_S\colon G\to S^{-1}$ be a morphism 
        of topological groups and
	$f\colon G\to \K[G]$  the coextension of $f_S$. Then by the
	universal property of $\K[G]$ there is a unique morphism 
	$f'\colon \K[G]\to S$ of weakly complete unital algebras
	such that $f'\circ \eta_G=f$, implying that 
        $(f'|S)\circ \eta_G^{\rm o}
	=f_S$ with the corestriction $\eta_G^{\rm o}\colon G\to S$
	of $\eta_G$ to $S$. Thus $S$ has the universal property
	of $\K[G]$; then the uniqueness of $\K[G]$ implies $S=\K[G]$.
\end{Proof}

\bsk

We recall from \cite{book3}, Corollary 2.29{\color{\red}(ii)} that it can
be proved in any theory of compact groups at a very early stage
that 

\ssk

\nin{\it every compact group has an isomorphic copy in the group of
units of a weakly complete unital algebra.}

\ssk

\nin As a consequence we have

\msk

\begin{Theorem} \label{th:gr-alg-comp-gr}
{\rm (The Group Algebra of a Compact Group)}
If $G$ is a compact group,
then $\eta_G\colon G\to\R[G]^{-1}$ induces an isomorphism of topological groups
onto its image.
\end{Theorem}

\nin In other words, {\it

\ssk

\nin any compact group may be considered as a subgroup of
the group of units of its weakly complete real group algebra.}

 \bsk

\subsection{The Group Algebra Functor $\K[-]$ is Multiplicative}
If $A$ and $B$ are weakly complete algebras, we have $(a_1\otimes b_1)(a_2\otimes b_2)
=a_1a_2\otimes b_1b_2$ which implies 
$$A^{-1}\otimes B^{-1}\subseteq (A\otimes B)^{-1},$$
where we have used the natural inclusion function
$j\colon A\times B\to A\otimes B$ and write $A^{-1}\otimes B^{-1}$
in place of $j(A^{-1}\times B^{-1})$.

\msk

Now let $G$ and $H$ be topological groups. Then 
$$\eta_G(G)\otimes \eta_H(H)\subseteq \K[G]^{-1}\otimes \K[H]^{-1}
    \subseteq \(\K[G]\otimes \K[H]\)^{-1},$$
and so we have the morphism 
$$G\times H\to \(\K[G]\otimes \K[H]\)^{-1},$$
$(g,h)\mapsto \eta_G(g)\otimes \eta_H(h)$ which, by the univeral property
of $\K[-]$ gives rise to a unique morphism 
$\alpha \colon \K[G\times H]\to \K[G]\otimes \K[H]$ such that
$$ (\forall (g,h)\in G\times H)\,
   \alpha(\eta_{G\times H}(g,h))=\eta_G(g)\otimes \eta_H(h).\leqno(1) $$
On the other hand, the morphisms $j_G\colon G\to G\times H$, 
$j_G(g)=(g,1_H)$ and $p_G\colon G\times H\to G$, $p_G(g,h)=g$ 
yield $p_Gj_G=\id_G$. Therefore 
$\K[p_G]\colon \K[G\times H]\to \K[G]$ is an algebra retraction,
and via $\K[j_G]$ we may identify $\K[G]$ with a subalgebra 
of $\K[G\times H]$; likewise $\K[H]$ is an algebra
retract of the algebra $\K[G\times H]$. Since
$(g,1)(1,h)=(g,h)$ in $G\times H$, with the identifications
of $\K[G], \K[H]\subseteq \K[G\times H]$ we have
$$ (\forall (g,h)\in G\times H)\, \eta_G(g)\eta_H(h)=
                \eta_{G\times H}(g,h) \in \K[G\times H].\leqno(2)$$
The function 
$$\K[G]\times \K[H]\to \K[G\times H],\quad 
(a,b)\mapsto ab$$
is a continuous bilinear map of weakly complete vector spaces;
therefore the universal property of the tensor product in 
$\cal W$ yields a unique $\cal W$-morphism

\cen{$\beta\colon \K[G]\otimes \K[H]\to \K[G\times H]$}

\noindent such that 
$$(\forall a\in \K[G],\,b\in \K[H]) \, \beta(a\otimes b)=
ab\in \K[G\times H].\leqno(3)$$
Now if for an arbitrary element $(g,h)\in G\times H$ we set 
$a=\eta_G(g)$ and $b=\eta_H(h)$,
then we have
$$\beta\(\eta_G(g)\otimes\eta_H(h)\)=a\otimes b=
ab=\eta_G(g)\eta_H(h)=\eta_{G\times H}(g,h).\leqno(4)$$
By Proposition \ref{p:generation}, $\eta_G(G)$ generates $\K[G]$ as 
weakly complete unital algebra and likewise
$\eta_H(H)$ generates $\K[H]$ in this fashion, and
the algebraic tensor product of $\K[G]$ and $\K[H]$
is dense in $\K[G]\otimes \K[H]$. 
Therefore, (4) implies 
$\beta\circ \alpha=\id_{\K[G\times H]}$. 
In other words, the diagram
$$
\begin{matrix}\K[G\times H]&\mapright{\id_{\K[G\times H]}}&\K[G\times H]\\
\lmapdown{\alpha}&&\mapup{\beta}\\
\K[G]\otimes \K[H]&\lmapright{\id_{\K[G]\otimes \K[H]}}&\K[G]\otimes \K[H]\\
\end{matrix}
$$ 

\medskip

\noindent
commutes.
Similarly, let us look at 
$\alpha\circ\beta:\K[G]\otimes \K[H]\to \K[G]\otimes \K[H]$: 
We recall (4) and (1) and verify
$$\alpha\(\beta(\eta_G(g)\otimes\eta_H(h)\)=\alpha(\eta_{G\times H}(g,h))
=\eta_G(g)\otimes\eta_H(h)$$
By the same argument as above we conclude 
$\alpha\circ\beta= \id_{\K[G]\otimes \K[H]}$.

Taking everything together, we have proved the following important result:
\msk
\begin{Theorem} \label{th:multiplicative} 
{\rm (Multiplicativity of the Group Algebra
Functor $\K[-]$)}
For two arbitrary topological groups $G$ and $H$
the natural morphisms of weakly complete unital algebras 
$\alpha\colon\K[G\times H]\to\K[G]\otimes\K[H]$ and
$\beta\colon\K[G]\otimes\K[H]\to\K[G\times H]$
are isomorphisms which are inverses of each other.
\end{Theorem}
 
\bsk

\subsection{Multiplication and Comultiplication on
the Group Algebra $\K[G]$} 
Let $G$ be a topological group and $\delta_G\colon G\to G\times G$
the diagonal morphism $\delta_G(g)=(g,g)$. Together with the constant
morphism $k_G\colon G\to \IdObject=\{1\}$ we have a comonoid $(\delta_G,k_G)$
 according to
the example following Definition \ref{d:comonoids}. Since the group-algebra functor
$\K[-]$ is multiplicative we have {\it morphisms of weakly complete
unital algebras} $\K[\delta_G]\colon \K[G]\to \K[G\times G]$
and $\K[k_G]\colon \K[G]\to \K[\{1\}]=\K$. By Theorem 
\ref{th:multiplicative} above
we have an isomorphism $\alpha_G\colon\K[G\times G]\to \K[G]\otimes\K[G]$
which gives us the following observation:

\msk

\begin{Lemma} For any topological group $G$, the weakly complete
group algebra $\K[G]$ supports a cocommutative and coassociative
comultiplation 
$$\gamma_G\colon \K[G]\to \K[G]\otimes \K[G],\quad
\gamma_G=\alpha_G\circ \K[\delta_G]$$ which is a morphism
of weakly complete unital algebras, and there is a
co-identity $\kappa_G\colon\K[G]\to\K$ which is an algebra morphism.
\end{Lemma}

\msk

The following fairly immediate remark will be relevant:

\begin{Remark} \label{r:grlike} 
Let $G$ be a topological group and $x\in\K[G]$.
If $x\in \eta_G(G)$, then the following statement holds:
$$\gamma_G(x)=x\otimes x\mbox{ and }\kappa(x)=1.\leqno(\dag)$$
The set of elements satisfying $(\dag)$ is linearly independent.
\end{Remark}

\msk

\begin{Proof} 
We recall the definition of 
$$\gamma_G=\left(\K[G]\mapright{\K[\delta_G]}\K[G\times G]\mapright{\alpha_G}\K[G]\otimes\K[G]\right).
\leqno(*)$$
If $a=\eta_G(g)$ for some $g\in G$, then $c(a)=\alpha_G(a,a)=a\otimes a$ by $(*)$
and by (1) above. 
The linear independence is a an exercise in linear algebra which one finds
in \cite{hofbig}, pp.~66, 67.
\end{Proof}

\msk

For each topological group $G$ the {\it opposite group} $G^{\rm op}$ is
the underlying topological space of $G$ together with the multiplication
$(g,h)\mapsto g{*}h$ defined by $g{*}h=hg$. The groups $G$ and $G^{\rm op}$
are isomorphic under the function $\inv_G\colon G\to G^{\rm op}$,
$\inv_G(g)=g^{-1}$. Analogously, every topological algebra $A$ gives
rise to an opposite algebra $A^{\rm op}$ on the same underlying topological
vectorspace but with the multiplication defined by $a{*}b=ba$, giving
us 
$$(A^{-1})^{\rm op}=(A^{\rm op})^{-1}$$
by definition, but not necessariy being isomorphic to $A$.
Consequently, $$((\K[G])^{-1})^{\rm op}=(\K[G]^{\rm op})^{-1}$$ and there are
 morphisms of topological groups $\eta_G\colon G\to \K[G]^{-1}$ 
and $\eta_{G^{\rm op}}\colon G^{\rm op}\to \K[G^{\rm op}]^{-1}$.
Accordingly, we have an isomorphism  
$\K[\inv_G]\colon\K[G]\to\K[G^{\rm op}]$ of weakly complete topological algebras
and, accordingly, an involutive isomorphism of topological groups
$\K[\inv_G]^{-1}\colon\K[G]^{-1}\to\K[G^{\rm op}]^{-1}$. 
This gives us a commutative diagram
$$
\begin{matrix}
 G&\mapright{\eta_G}&\K[G]^{-1}\\
\lmapdown{\inv_G}&&\mapdown{\K[\inv_G]^{-1}}\\
G^{\rm op}&\lmapright{\eta_{{G^{\rm op}}}}&\K[G^{\rm op}]^{-1}.\\
\end{matrix}
$$
producing an isomorphism of weakly complete algebras
$\K[G]\to \K[G^{\rm op}]$.

But we also have a commutative diagram

$$
\begin{matrix}G&\mapright{\eta_G}&\K[G]^{-1}\hfill\\
\lmapdown{\inv_G}&&\mapdown{\inv_{\K[G]^{-1}}}\hfill\\
G^{\rm op}&\lmapright{\eta_{{G^{\rm op}}}}&(\K[G]^{-1})^{\rm op}
     =(\K[G]^{\rm op})^{-1}.\\
\end{matrix}
$$
Let us abbreviate 
$f\defi\(\eta_{G^{\rm op}})\circ\inv_G\colon G\to(\K[G]^{\rm op})^{-1}$.
So by the adjunction formalism, there is a unique involutive
isomorphism
$f'\colon \K[G]\to \K[G]^{\rm op}$ of weakly complete
algebras such that
$f=f'|\K[G]^{-1}\circ\eta_G$. 

\bsk

We have a  grounding functor $A\to|A|$ from the cagegory
$\WA$ of weakly complete algebras to the category $\W$ of weakly
complete vector spaces, where $|A|$ is simply the weakly complete
vector space underlying the weakly complete algebra $A$.
With this convention we formulate the following definition:

\begin{Definition} \label{l:symmetry} For each topological
group $G$ there is a morphism of weakly complete vector spaces
$$\sigma_G\defi |f'|\colon |\K[G]| 
            \to |\K[G]^{\rm op}|=|\K[G]|, $$
often called {\it symmetry} {\color{\red} or {\it antipode}}
such that for any topological group $G$ we have
$$(\forall g\in G)\, \sigma_G(\eta_G(g))=\eta_G(g^{-1})=\eta_G(g)^{-1}.$$ 
\end{Definition}

\nin Equivalently, 
$$(\forall x\in \eta_G(G))\, \sigma_G(x)\.x=x\.\sigma_G(x)=1,$$
and, using the bilinearity of multiplication
$(x,y)\mapsto xy$ in $\K[G]$ and
 defining $\mu_G\colon \K[G]\otimes \K[G]\to \K[G]$ by $\mu_G(x\otimes y)=xy$
and remembering from Remark \ref{r:grlike} that
$\gamma_G(x)=x\otimes x$ for all $x=\eta_G(g)$ {\color{\red} with}
 some $g\in G$, once more, equivalently,
$$(\forall x \in\eta_G(G))\, 
(\mu_G\circ(\sigma_G\otimes \id_{\K[G]})\circ\gamma_G)(x)=1.\leqno(+)$$
By Proposition \ref{p:generation}, the weakly complete algebra $\K[G]$
is the closed linear span of $\eta_G(G)$, equation $(+)$ holds
in fact for all elements or $\K[G]$.
Thus we have shown the following 

\msk

\begin{Proposition} \label{p:prelim-hopf} For any topological group $G$, 
the following diagram involving natural morphisms of weakly 
complete vector spaces commutes where $W_G=|\K[G]|$:
$$
\begin{matrix}
W_G\otimes W_G&\hfill\mapright{\sigma_G\otimes \id_{W_G}}\hfill&W_G\otimes W_G\\
\lmapup{\gamma_G}&&\mapdown{\mu_G}\\
W_G&\lmapright{\iota_G \circ |\kappa_G|}&W_G,\\
\end{matrix}
$$
\msk \noindent
where $\kappa_G\colon \K[G]\to \K=\K[\IdObject]$, $\IdObject=\{1\}$  is 
induced by the constant morphism and $\iota_G\colon \K=\K[\IdObject]\to\K[G]$ 
by the unique inclusion $\IdObject\to G$.

In this diagram, the multiplication $\mu_G$ defines the algebra structure
of the group algebra $\K[G]$ on $W_G$ and the comultiplication
$\gamma_G$ is an algebra morphism.
\end{Proposition}

\msk

\subsection{Weakly Complete Bialgebras and Hopf-Algebras}
We need to view the basic result of Proposition \ref{p:prelim-hopf} in
its more general abstract frame work.

\nin
In an arbitrary monoidal category $({\cal A},\otimes)$ we have
monoids $A\otimes A\mapright{m} A\mapleft{u} \IdObject$ and
comonoids $\IdObject\mapleft{k} C\mapright{c} C\otimes C$.
Every pair of monoids $A$ and $B$ gives rise to a
monoid $A\otimes B$. In particular, for each monoid $A$
also $A\otimes A$ is a monoid. An $\cal A$-morphism
$f\colon A\to B$ between monoids is a {\it monoid
morphism} if they respect the monoid structure
in an evident fashion. These matters are discussed
in a comprehensive fashion in \cite{book3} in Appendix 3
in the section entitled ``Commutative Monoidal Categories
and their Monoid Objects''; see pp.787 ff., Example A3.61ff.

\msk

\begin{Definition} \label{d:hopf}
 (a) A {\it bimonoid} in a commutative
monoidal category is an object together with both a monoid
structure $(m,u)$ and comonoid structure $(c,k)$,
$$ A\mapright{c}A\otimes A\mapright{m}A\hbox{\quad and\quad}
 \IdObject\mapright{{\color{\red}u}}A\mapright{{\color{\red}k}}\IdObject,$$
such that $c$ is a monoid morphism.

\msk

\nin
(b) A {\it group} (or often {\it group object} in
a commutative monoidal category) is a bimonoid with commutative
comultiplication and with an $\cal A$-morphism
$\sigma\colon A\to A$, called {\it inversion} or
{\it symmetry} (as the case may be) which
makes the following diagram commutative
$$
\begin{matrix}
A\otimes A&\mapright{\sigma\otimes\id}&A\otimes A\\
\lmapup{c}&&\mapdown{m}\\
A&\lmapright{uk}&A,\\
\end{matrix}\leqno(\sigma)
$$
plus a diagram showing its compatibility
with the comultiplication (see \cite{book3},
Definition A3.64.ii).

{\color{\red}

\nin
(c) In our commutative monoidal categories $({\mathcal V},\otimes)$
and $({\mathcal W},\otimes)$ of $\K$-vector spaces, respectively,
weakly complete $\K$-vectorspaces, a group object $(A,m,c,u,k,\sigma)$
is called a {\it Hopf algebra}, respectively, {\it a weakly complete
Hopf algebra}.}   
\end{Definition}

\msk
In reality, the definition of a bimonoid is symmetric and the equivalent
conditions that $c$ be a monoid morphism, respectively, that
$m$ be a comonoid morphism can be expressed in one
commutative diagram (see \cite{book3}, Diagram following
Definition A3.64, p.~793). Also it can be shown that in a group 
the diagram arising from the diagram $(\sigma)$ by replacing
$\sigma\otimes \id$ by $\id\otimes\sigma$ commutes as well.

\msk

{\color{\red}
Since the group objects and group morphisms in a commutative monoida;
category form a category (see e.g.\ \cite{book3}, paragraph
following Exercise EA3.34, p.~793), the Hopf-algebras in
$({\mathcal V},\otimes)$ and $({\mathcal W},\otimes)$ do constitute
categories in their respective contexts.}  

\msk

\begin{Example} (a) In the category of sets (or indeed
in any category with finite products and terminal object)
every monoid is automatically a bimonoid if one takes
as comultiplication the diagonal morphism $X\to X\times X$.

\msk

(b) A bimonoid in the category of 
vector spaces with the tensor product as multiplication
 is called a {\it bialgebra}. Sometimes it is also
called a {\it Hopf-algebra}, and sometimes this
name is reserved for a group object in this category.
The identity $u\colon \K\to A$ of an algebra identifies
$\K$ with the subalgebra $\K\.1$ of $A$.
The coidentity $k\colon A\to\K$ of a bialgebra {\color{\red} may
be} considered as the endomorphism $uk\colon A\to A$
mapping $A$ onto the subalgebra $\K\.1$ and then is
often called the {\it augmentation}.
\end{Example}

In the present context we do restrict our attention to the 
commutative monoidal categories
$({\cal V},\otimes)$ and $({\cal W},\otimes)$ of
$\K$-vector spaces, respectively, weakly complete
$\K$-vector spaces with their respective tensor products.
In terms of the general terminology which is now available 
to us, the result in Proposition \ref{p:prelim-hopf} 
may be rephrased as follows:

\bsk

\begin{Theorem} \label{th:hopf} 
{\rm (The Hopf-Algebra Theorem for Weakly Complete
Group Algebras)} For every topological group $G$, the group
algebra $\K[G]$ is a weakly complete Hopf-algebra with a comultiplication $\gamma_G$
and symmetry $\sigma_G$.
\end{Theorem} 

\bsk

\bsk
\section{Grouplike and Primitive Elements}
\bsk

In any theory of Hopf algebras it is common to single out
 two types of special elements, and
we review them in the case of weakly complete Hopf algebras.

\msk

\begin{Definition} \label{d:grouplike} Let $A$ be a 
weakly complete coassociative
coalgebra with comultiplication $c$ and coidentity $k$. Then an element 
$a\in A$ is called {\it grouplike} if $k(a)=1$ and
$c(a)=a\otimes a$. 

\msk

If $A$ is a bialgebra, $a\in A$ is called {\it primitive}, if
$c(a)= a\otimes 1 + 1\otimes a$.
\end{Definition}

For any $a\in A$ with $c(a)=a\otimes a$, the conditions $a\ne0$ 
and $k(a)=1$ are equivalent.

\msk

These definitions apply, in particular, to any weakly complete
Hopf algebra and thus especially to each weakly complete
group algebra $\K[G]$. By Remark \ref{r:grlike} the
subset $\eta_G(G)$ is a linearly independent set  of grouplike
elements.

\begin{Remark} \label{r:warning} In at least one source 
on bialgebras in earlier
contexts, the terminology conflicts with the one introduced here
which is now commonly accepted. In \cite{hofbig}, p.~66, Definition
10.17, the author calls a grouplike element  in a coalgebra
{\it primitive}. Thus some caution is in order concerning
terminology. Primitive elements in the
sense of Definition \ref{d:grouplike} do not occur in \cite{hofbig}.
\end{Remark}

\msk
The following observations {\color{\red} are proved
 straightforwardly:}

\begin{Lemma} \label{l:monoid} The {\color{\red} set $G$ 
of}  grouplike elements of a weakly complete
bialgebra $A$ is a closed submonoid of $(A,\.)$ 
and {\color{\red}  the set $L$ of}  primitive elements
of $A$ is a closed Lie subalgebra of $A_{\rm Lie}$.
{\color{\red} If $A$ is a Hopf algebra, then $G$ is a closed subgroup
of $A^{-1}$.}
\end{Lemma}

\bsk

For a morphism $f\colon W_1\to W_2$ of weakly complete
vector spaces let 
$f'=\Hom(-,\K)\colon W_2'\to W_1'$ denote the dual 
morphism of vector spaces.

For a weakly complete coalgebra $A$ let $A'
= \Hom(A,\K)$ be the dual of $A$. 
Then $A'$ is an algebra: If 
$c\colon A\to A\otimes A$ is its
comultiplication,
then $c'\colon A'\otimes A'\to A'$ is the
 multiplication of $A'$.
For a unital algebra $R$ and a weakly complete coalgebra $A$
in duality let $(a,g)\mapsto\<a,g\>:R\times A\to\K$ 
denote the pairing between $R$ and $A$, where for $f\in\R=\Hom(A,\K)$ and
$a\in A$ we write $\<f,a\>=f(a)$.

{\color{\red} 
\begin{Definition} \label{d:sixfour} Let $R$ be a unital algebra over $\K$. Then
a {\it character} of $R$ is a morphism of unital algebras
$R\to \K$. The subset of $\K^R$ consisting of all algebra morphisms
 inherits the topology of 
pointwise convergence from $\K^R$ and as a topological space
is called the {\it spectrum}
of $R$ and is denoted $\Spec(R)$.

If $k\colon R\to \K$ is a morphism of algebras, then 
an element $d\in{\cal V}(R,\K)$ is called a \emph{derivative} (sometimes also 
called a \emph{derivation} or  \emph{infinitesimal character}) of $R$ (with
respect to $k$) if it satisfies 
$$ (\forall x, y\in R)\, d(xy)= d(x)k(y)+k(x)d(y).$$
The set of all derivatives of $R$ is denoted $\Der(R)$. 
\end{Definition}}
\msk 

 Now let $R$ be a unital algebra and $A\defi R^*$
its dual weakly complete coalgebra with comultiplication $c$ such 
that $ab=c'(a\otimes b)$ for all $a,b\in R$. In these circumstances
we have:
 
\begin{Proposition} \label{p:grouplike} 
Let $g\in A$. Then the following statements are equivalent:
\begin{enumerate}[\rm(i)]
	\item $g\in A$ is grouplike in the coalgebra $A$.
	\item $g\colon R\to \K$ is a character of $R$, that is,
               is an element of $\Spec(R)$.
\end{enumerate}
\end{Proposition}

\begin{Proof} 
The dual of $A\otimes A$ is
$R\otimes R$ in a canonical fashion such that for $r_1,r_2\in R$ and
$h_1,h_2\in A$ we have 
$$\<r_1\otimes r_2,h_1\otimes h_2\>=\<r_1,h_1\>\<r_2,h_2\>.$$
The set of linear combinations $L=\sum_{j=1}^n a_j\otimes b_j\in A\otimes A$ 
is dense in $A\otimes A$. So 
two elements $x,y\in R\otimes R$ agree if and only if for all such 
linear combinations $L$ we have 
$$\<x,L\>=\<y,L\>,$$
and this clearly holds if and only if for all $a,b\in A$ we have
$$\<x,a\otimes b\>=\<y,a\otimes b\>.$$
We apply this to $x=c(g)$ and $y=g\otimes g$
and observe that (i) holds if and only if 
$$(\forall r,s\in R)\,\<r\otimes s,c(g)\>=\<r\otimes s, g\otimes g\>.
\leqno(0)$$
Now
$$g(rs)=\<rs,g\>=\<m(r\otimes s),g\>=\<r\otimes s,c(g)\>.\leqno(1)$$
$$g(r)g(s)=\<r,g\>\<s,g\>=\<r\otimes s,g\otimes g\>.\leqno(2)$$
So in view of (0),(1) and (2), assertion (i) holds if and only if
$g(rs)=g(r)g(s)$ for all $r, s\in R$ is true. Since $g$ is a 
linear form on $A$, this means exactly that a nonzero $g$ is 
a morphism of weakly complete algebras, i.e., 
$g\in\Spec(R)$.
\end{Proof}

\bsk
A particularly relevant consequence of Proposition \ref{p:grouplike}
will be the following con{\color{\red}clusion} for the group algebra $\R[G]$ of a
{\it compact} group. Recall that after Theorem \ref{th:gr-alg-comp-gr}
we always identify a compact group $G$ with a subgroup of $\R[G]$
in view of the embedding $\eta_G$. We shall see that the set $G\subseteq\R[G]$
contains already {\it all} grouplike elements (see Theorem 
\ref{th:comp-group-like} below).

\msk

Let us return  to the primitive elements of a unital bialgebra. {\color{\red}
For this purpose assume that $R$ is not only a unital algebra, but a bialgebra
over $\K$, which  implies  that its dual $A\defi R^*$ is not only a coalgebra but a bialgebra
as well.}
 
\begin{Proposition} \label{l:primitive} Let $R$ be a unital bialgebra 
and  $d\in A$. 
Then the following statements are equivalent:
\begin{enumerate}[\rm(i)]
	\item $d$ is primitive in the bialgebra $A$.
	\item $d\colon R\to\K$ a derivative of $R$ with respect to the coidentity $k$, \
              that  is, an element in $\Der(R)$.
\end{enumerate}
\end{Proposition}

The procedure of the proof of Proposition \ref{p:grouplike} 
allows us to leave the explicit proof of this proposition 
as an exercise.

\msk

\begin{Definition} \label{d:grlike}  Let $A$ be a weakly complete Hopf algebra.
Then we write $\Gamma(A)$ for the subset of all grouplike
elements of $A$ and $\Pi(A)$ for the subset of all primitive elements.
\end{Definition}

In view of Proposition \ref{p:monoids}, 
the following proposition is now readily verified:

\begin{Proposition} \label{l:substructures} 
In a weakly complete Hopf algebra $A$, the subset $\Gamma(A)$ is
a closed subgroup of $A^{-1}$ and $\Pi(A)$ is 
 a closed Lie subalgebra of $A_{\rm Lie}$. If $f\colon A\to B$
is a morphism of weakly complete  Hopf algebras, then
 $f(\Gamma(A))
\subseteq \Gamma(B)$, and $f(\Pi(A))\subseteq \Pi(B)$.  
\end{Proposition} 

\msk

{\color{\red}
Accordingly, if $\mathcal{WH}$ denotes the category of weakly complete
Hopf algebras, then $\Gamma$ defines a functor from $\mathcal{WH}$ to
the category of topological groups $\mathcal{T\hskip-1pt G}$, also called
{\it the grouplike grounding}, while $\Pi$ defines 
a functor from $\mathcal{WH}$ to the category of topological 
Lie algebras over $\K$.

The functor $\HH\defi(G\mapsto \R[G]): \mathcal{T\hskip-1pt G}\to \mathcal{WH}$ 
from topological groups to weakly complete Hopf algebras was known to us
since Theorem \ref{th:wcga-thm} as a functor into the bigger category of all 
weakly complete unital algebras. That theorem will now be sharpened
as follows:

\begin{Theorem} \label{p:grounding} {\rm (The Weakly Complete Group Hopf
Algebra Adjunction Theorem)}\quad 
 The functor $\HH\colon \mathcal{T\hskip-1pt G}\to \mathcal{WH}$
from topological groups to weakly complete Hopf algebras is left adjoint
to the grouplike grounding $\Gamma\colon \mathcal{WH}\to \mathcal{T\hskip-1pt G}$.
\end{Theorem} 

In other words, for a topological group $G$ there is natural
morphism of topological groups $\eta_G\colon G\to \Gamma(\HH(G))=
\Gamma(\K[G])$ such that for each morphism of topological groups
$f\colon G\to \Gamma(A)$ for a weakly complete Hopf algebra $A$  
there is a unique morphism of weakly complete Hopf algebras 
${\color{\red}f'}\colon \HH(G)\to A$ such that $f(g)=f'(\eta_G(g))$ for all
$g\in G$:

$$ 
\begin{matrix}&\mathcal{T\hskip-1pt G}&&\hbox to 7mm{} &\mathcal{WH}\\
\noalign{\vskip3pt}
\noalign{\hrule}\\
\noalign{\vskip3pt}%
        \mdvertspace G&\mapright{\eta_G}&\Gamma(\HH(G))&\hbox to 7mm{} &
               \HH(G)=\K[G]\\
\lmapdown{\forall f}&&\mapdown{\Gamma(f')}&\hbox to 7mm{}&
            \mapdown{\exists! f'}\\
\muvertspace \Gamma(A)&\lmapright{\id}&\Gamma(A)&\hbox to 7mm{}&A\\
\end{matrix}
$$
}

\msk

\begin{Proof} Let $A$ be a weakly complete Hopf algebra and
$f\colon G\to \Gamma(A)$ a continuous group morphism. Since
$A$ is in particular a weakly complete associative unital algebra and
$\Gamma(A)\subseteq A^{-1}$, by the Weakly Complete Group Algebra Theorem
\ref{th:wcga-thm} there is a unique morphism $f'\colon\K[G]\to A$
of weakly complete algebras such that $f(g)=f'(\eta_G(g))$
for all $g\in G$. Since each $\eta_G(g)$ is grouplike by
Remark \ref{r:grlike} we have $\eta_G(G)\subseteq \Gamma(\HH(G))$.
We shall see below that the morphism $f'$ of 
weakly complete algebras is indeed a morphism of weakly complete
 Hopf algebras and therefore maps grouplike element into
grouplike elements. Hence
 $f'$ maps $\Gamma\(\HH(G)\)$ into $\Gamma(A)$. 

We now  have to show that $f'$ is a morphism of Hopf algebras, that is,

(a)  it respects comultiplication,

(b)  coidentity, and

(c)  symmetry.

\nin
For (a) we must show that the following diagram commutes:

$$\begin{matrix} \K[G]&\mapright{c_{\K[G]}}&\K[G]\otimes\K[G]&\cong&
                                              \K[G\times G]\cr
                  \lmapdown{f'}&&\mapdown{f'\otimes f'}&&\cr
                   A&\mapright{c_A}&A\otimes A. &&\end{matrix}$$
Since $\K[G]$ is generated as a topological algebra by $\eta_G(G)$
by Proposition \ref{p:generation}, it suffices to track
all elements $x=\eta_G(g)\in \K[G]$ for $g\in G$. Every
such element is grouplike in $\K[G]$ by Remark \ref{r:grlike}, and so 
$(f'\otimes f')c_{\K[G]}(x)= (f'\otimes f')(x\otimes x)=f'(x)\otimes f'(x)$
in $A\otimes A$, while on the other hand $f'(x)=f(g)\in
\Gamma(A)$, whence $c_A(f'(x))=f'(x)\otimes f'(x)$ as well. 
This proves (a).

For (b) we must show that the following diagram commutes:
$$\begin{matrix} \K[G]&\mapright{k_{\K[G]}}& \K&\cong&
                                              \K[\{1_G\}]\cr
                  \lmapdown{f'}&&\mapdown{\id_\K}\cr
                   A&\mapright{k_A}&\K. &&\end{matrix}$$
Again it suffices to check the elements $x=\eta_A(g)$.
Since all grouplike elements are always mapped to 1,
this is a trivial exercise.

Finally consider (c), where again we follow all elements
$x=\eta_A(g)$. On the one hand we have $f'(\sigma_{\K[G]}(x))
=f'(x^{-1})=f'(x)^{-1}$ in $A^{-1}$. But again $f'(x)$
is grouplike, and thus $\sigma_A(f'(x))=f'(x)^{-1}$, which
takes care of case (c), and this completes the proof of
the theorem.
\end{Proof}

\msk

As with any adjoint pair of functors there is an alternative
way to express the adjunction in the preceding theorem:
see e.g.\ \cite{book3}, Proposition A3.36, p.~777:

\begin{Corollary} \label{c:grounding} For each weakly complete 
Hopf algebra $A$ there  is a natural morphism {\color{\red}
of Hopf algebras}
$\epsilon_A\colon \HH(\Gamma(A))\to A$ such that for
any topological group $G$ and any morphism of {\color{\red}
Hopf algebras}
$\phi\colon \HH(G)\to A$
there is a unique continuous group morphism 
$\phi'\colon G\to \Gamma(A)$ such that for each $x\in \HH(G)=\K[G]$
one has $\phi(x)=\epsilon_A(\K[\phi](x))$, where $\K[\phi]=\HH(\phi)$.
\end{Corollary}

 The formalism we described in Corollary \ref{c:alternate} we can
formulate with the present adjunction as well:
                                                                   
\begin{Corollary} \label{c:alternate2}
For any weakly complete Hopf algebra
 $A$ and any topological group $G$ we have

$$ (\forall A)\, \Big(\Gamma(A)\quad\mapright{\eta_{\Gamma(A)}}\quad 
       \Gamma(\K[\Gamma(A)])
 \quad \mapright{\Gamma(\epsilon_A)}\ \Gamma(A)\Big)
    =\id_{\Gamma(A)},\hbox{ and}\leqno(1)$$

\medskip

$$(\forall G)\,\Big(\K[G]\quad\mapright{\K[\eta_G]}\quad 
       \K\big[\Gamma(\K[G])\big]
\quad \mapright{\epsilon_{\K[G]}}\quad \K[G]\Big)=\id_{\K[G]}.\leqno(2)$$
\end{Corollary}

{\color{\red}

\begin{Definition}\label{d:span-grl} For any weakly complete Hopf 
algebra $A$ we let $\SS(A)$ denote the closed linear
span $\overline{\Span}\(\Gamma(A)\)$ in $A$.
\end{Definition}
Since $\Gamma(A)$ is a multiplicative subgroup of $(A,\cdot)$,
clearly $\SS(A)$ is a weakly closed subalgebra of $A$.
Similarly, one easily verifies $c\(\Gamma(A)\)\subseteq \Gamma(A\otimes A)
\subseteq \SS(A\otimes A)$ for the comultiplication of $A$,
and so $c\(\SS(A)\) \subseteq \SS(A\otimes A)$. Therefore

\ssk

\nin
{\it $\SS(A)$ is  Hopf subalgebra of $A$.}

\msk
\begin{Lemma} \label{l:generation2} If $A=\K[G]$ is the weakly complete group algebra
of a topological group, then $A=\SS(A)$.
\end{Lemma}

\begin{Proof} By Remark \ref{r:grlike},
$$\eta_G(G)\subseteq \Gamma(A).\leqno(1)$$
By Proposition \ref{p:generation},
$$\overline{\Span}\(\eta_G(G)\)=A.\leqno(2)$$
Hence  
$$\SS(A)=\overline{\Span}\(\Gamma(A)\)=A.\leqno(3)$$ 
\vglue-17pt
\end{Proof}
 
\begin{Corollary} \label{c:epsilon-epic} In the circumstances of
{\rm Corollary \ref{c:grounding}}, 
$$\im(\epsilon_A)=\epsilon_A\(\HH\(\Gamma(A)\)\)=\SS(A).$$
\end{Corollary}

\begin{Proof} Set $B=\HH(\Gamma(A))$; then $\epsilon_A\colon B\to A$
is a morphism of Hopf algebras, and so, in particular, a morphism  of
weakly complete vector spaces. Hence $\im(\epsilon_A)=
\epsilon_A(B)$ is a closed Hopf subalgebra of $A$. Since
$\epsilon_A$ is a morphim of Hopf algebras,
$\epsilon_A(\Gamma(B))\subseteq \Gamma(A)$ and thus
$\epsilon_A\(\SS(B)\)\subseteq \SS(A)$. By Lemma \ref{l:generation2} 
we have $B=\SS(B)$ and so $\im(\epsilon_A)=\epsilon_A\(\SS(B)\)$
and so  $\im(\epsilon_A) \subseteq\SS(A)$.

On the other hand, by Corollary \ref{c:alternate2},
we have $\Gamma(A)\subseteq \im(\epsilon_A)$ and since
$\im(\epsilon_A)$ is closed, we conclude
$\SS(A)\subseteq \im(\epsilon_A)$ which completes the proof.
\end{Proof}

In particular, 
{\it $\epsilon_A$ is quotient homomorphism if and only if
$A=\SS(A)$.}

\msk

After we identified the image of the Hopf algebra morphism $\epsilon_A$, 
the question
arises for its kernel.
Thus let $J=\ker\epsilon_A$. As a
kernel of a morphism of Hopf algebras, $J$ is, firstly, the
kernal of a morphism of weakly complete algebras, and
secondly,  satisfies the
condition
$$  c(J)\subseteq A\otimes J + J\otimes A.$$
(See e.g.\ \cite{hofbig}, pp.\ 50ff.,  notably, Proposition 10.6
with appropriate modifications.)
We leave the significance of this kernel open for the time being
and return to it in the context of compact groups.
}

Here we record what is known in general on the pro-Lie group
$\Gamma(A)\subseteq A^{-1}$ for a weakly complete Hopf algebra $A$
in general.

\bsk

It was first observed in \cite{dah} that any weakly complete 
Hopf algebra $A$ accommodates grouplike elements (and so the spectrum
$\Spec(R)$ of the dual Hopf algebra of $A$) and its primitive 
elements (and hence the derivatives $\Der(R)$ of the dual) 
in a particularly satisfactory set-up for whose proof we refer
to \cite{dah}, \cite{dahii}, or \cite{hofmor}. In conjunction
with the present results, we can formulate the following theorem:

\begin{Theorem} \label{th:hopf-prolie}
Let $A$ be a weakly complete $\K$-Hopf-algebra and $A^{-1}$
its group of units which is an almost connected pro-Lie group 
{\rm (Theorem \ref{th:fourth})} {\color{\red}
containing as a closed subgroup  the group $\Gamma(A)$ of
grouplike elements. In particular, $\Gamma(A)$ is a pro-Lie group. 
Correspondingly, the set $\Pi(A)$ of primitive elements of $A$ is
a closed
Lie subalgebra of the pro-Lie algebra 
$A_{\rm Lie}$ {\rm(Theorem \ref{th:second})} and may be identified
with the Lie algebra $\L(\Gamma(A))$ of $\Gamma(A)$. The exponential function
$\exp_A\colon A_{\rm Lie}\to A^{-1}$ restricts to the exponential
function $\exp_{\Gamma(A)}:\L(\Gamma(A))\to \Gamma(A)$.} 
\end{Theorem}

Again, we shall pursue this line {\color{\red} below} by applying it specifically  to
the  real group Hopf algebras $A\defi \R[G]$ in the special
case that $G$ is a {\it compact} group by linking it with
particular information we have on compact groups.
But first we have to explore duality more explicitly on the 
Hopf algebra level.

{\color{\red}

       \section {The Dual  of  Weakly Complete Hopf Algebras
                        and of Group Algebras}

}

In this section we will have a closer look at the dual space $A'$ of a weakly 
complete {\color{\red} Hopf} algebra $A$. We let $G$ denote the topological
group  of grouplike elements $g\in A$.
The underlying weakly complete vector space of 
$A$ is a topological left and
right $G$-module $A$ with the module operations
$$
\begin{matrix}
&(g,a)\mapsto g\.a: &G\times A\to A,\ &g\.a:= ga, &\mbox{and}\\
&(a,g)\mapsto a\.g: &G\times A\to A,\ &a\.g:= ag.\\
\end{matrix}
$$
Recall that $\I(A)$ is the filterbasis of closed two-sided ideals
$J$ of $A$ such that $A/J$ is a finite dimensional algebra
and that $A\cong\lim_{J\in\I(A)}A/J$.
We can clearly reformulate Corollary \ref{c:ideal-filter}
in terms of $G$-modules as follows:

\begin{lemma} \label{l:cofinite} 
	For the topological group $G=\Gamma(A)$,
	the $G$-module $A$ has a filter basis $\I(A)$ of closed two-sided submodules
	$J{\subseteq}A$ such that $\dim(A/J){<}\infty$ and that 
	$A=\lim_{J\in\I(A)} A/J$ is a strict projective limit of finite dimensional
        $G$-modules. The filter basis $\I(A)$ in $A$
	converges to $0\in A$.
\end{lemma}

For a $J\in\I(A)$ let 
$J^\perp=\{f\in A': (\forall a\in J)\,\<f,a\>=0\}$ 
denote
the annihilator of $J$ in the dual $V$ of $A$. We compare the ``Annihilator Mechanism''
from \cite{book3}, Proposition 7.62 and observe the following configuration:
$$
\begin{matrix}\mdddvertspace A&\hskip 2truecm&\{0\}&&\\               
        \Big|&&\Big|&\Bigg\}&\!\!\cong\ (A/J)'  \\                                                                \muuuvertspace J&&J^\perp&&\\                                       
         \Big|&&\Big|&\Bigg\}&\!\!\cong\ J'\\                              
          \{0\}&&\hphantom{.}A'.&&\\
\end{matrix}                                      
$$
\medskip

\noindent
In particular we recall the fact that $J^\perp\cong (A/J)'$
showing that $J^\perp$ is a finite-dimensional $G$-module 
on either side. By simply dualizing Lemma \ref{l:cofinite}
we obtain 

\begin{Lemma} \label{l:finite} 
For the topological group $G=\Gamma(A)$, the dual $G$-module 
$R\defi A'$ of the weakly complete $G$-module $A$ has an up-directed 
set $\D(R)$ of finite-dimensional two-sided $G$-submodules (and
$\K$-coalgebras!)  $F\subseteq R$ 
such that $R$ is the direct limit
\[
 R= \colim\limits_{F\in\D(R)}F=\bigcup_{F\in \D(R)}F.
\]
The colimit is taken in the category of (abstract) $G$-modules, 
i.e.~modules without any topology.
\end{Lemma}

\noindent This means that for the topological group $G=\Gamma(A)$, 
every element $\omega$
of the dual of $A'$ is contained in a finite dimensional 
left- and right-$G$-module (and $\K$-subcoalgebra).

We record this in the following form:

\begin{Lemma} \label{l:translates} 
 Let $\omega\in A'$. Then the vector subspaces $\Span(G\.\omega)$ and
$\Span(\omega\.G)$ of both the left orbit and the right orbit of 
$\omega$ are finite dimensional,
and both are contained in a finite dimensional $\K$-subcoalgebra of $A'$.
\end{Lemma}

\noindent For any $\omega\in A'$ the restriction $f\defi \omega|G:G\to \K$ 
is a continuos function such that each of the sets of translates 
$f_g$, $f_g(h)=f(gh)$,
respectively, ${}_gf$, ${}_gf(h)=f(hg)$ forms a finite 
dimensional vector subspace of the space $C(G,\K)$ of the vector space of all 
continuous $\K$-valued functions $f$ on $G$.

\msk
\begin{Definition} \label{d:R(G,K)}
For an arbitrary topological group $G$ we define
$\RR(G,\K)
\subseteq C(G,\K)$ to be that set of continuous functions 
$f\colon G\to \K$ for which the linear span
of the set of translations ${}_gf$, ${}_gf(h)=f(hg)$,
 is a finite dimensional vector subspace of $C(G,\K)$.
The functions in $\RR(G,\K)$ are called {\it representative
functions.}
\end{Definition}

In Lemma \ref{l:translates} we saw that for a weakly complete
Hopf algebra $A$  and its dual $A'$ (consisting of continuous
linear forms)  we have a natural linear
map 
$$\tau_A\colon A'\to \RR(\Gamma(A),\K),\ \tau_A(\omega)(g)=
(\omega|\Gamma(A))(g).$$
An element $\omega\in A'$ is in the kernel of $\tau_A$ if and only if
$\omega(\Gamma(A))=\{0\}$ if and only if $\omega\(\SS(A)\)=\{0\}$ if and only if
$\omega\in \SS(A)^\perp$.
We therefore observe:

\begin{Lemma} \label{l:tau-a} There is an exact sequence 
of $\K$-vector spaces
$$0\to\SS(A)^\perp\mapright{\inc}A'\mapright{\tau_A}\RR(\Gamma(A),\K).$$
\end{Lemma} 
Let us complete this exact sequence by determining the surjectivity
of $\tau_A$ at least
in the case that $A=\K[G]$ for a topological group $G$.

Indeed, let us show that a function $f\in \RR(G,\K)$ 
for an arbitrary 
topological group $G$ is in the image of $\tau_A$.

So let $f\in\RR(G,\K)$,  then
$U\defi \Span\{{}_gf:g\in G\}$ is an $n$-dimensional $\K$-vector subspace
of $C(G,\K)$. Accordingly,
 $M_U\defi\Hom(U,U)$ is an $n^2$-dimensional unital $\K$-algebra such that 
$\pi\colon G\to {\rm Gl}(U)$, $\pi(g)(u)={}_gu$, $u\in U$,
 is a continuous representation of $G$ on $U$.
Now we employ the universal property of the group algebra 
$\K[G]$ expressed
in Theorem \ref{th:wcga-thm} and find a unique morphism of
weakly complete algebras
$\pi'\colon \K[G]\to M_U$  inducing 
 a morphism of
topological groups $\pi'|\K[G]^{-1}:\K[G]^{-1}\to {\rm Gl}(U)$ 
such that $\pi=\(\pi'|\K[G]^{-1}\)\circ\eta_G$. If $\eta_G$ happens to
be an embedding and $G$ is considered as a subgroup of 
$\K[G]^{-1}$, then the representation $\pi$ is just the
restriction of $\pi'$ to $G$. Now
we recall that $f\in U$ and note that for any endomorphism
$\alpha\in M_U$ of $U$ we have
$\alpha(f)\in U\subseteq C(G,\K)$. Then for the identity $1\in G$
the element $\alpha(f)(1)\in\K$ is well defined and so
we obtain a linear map 
$L\colon M_U\to\K$ by defining
$L(\alpha)=\alpha(f)(1)$. So we finally define the linear form
$\omega_\pi=L\circ \pi'\colon \K[G]\to \K$ and calculate 
for any $g\in G$ the field element
$$\omega_\pi(\eta_G(g))=L\(\pi'(\eta_G(g))\)=
L(\pi(g))=\pi(g)(f)(1)={}_gf(1)=f(g).$$
Thus we have in fact shown the following 

\begin{Lemma} \label{l:rafael} For any topological group $G$ the following
statements are equivalent for a continuous function $f\in C(G,\K)$:
\begin{enumerate}[\rm(i)]
\item  $f\in\RR(G,\K)$.
\item There is a continuous linear form $\omega\colon \K[G]\to \K$
     such that $\omega\circ\eta_G=f$
\end{enumerate}
\end{Lemma}

Remember here that in case $G$ is naturally embedded into $\K[G]$ that
condition (ii) just says

\msk

\hskip-30pt(ii$'$) $\omega|G=f$ {\it for some} $\omega\in \K[G]'$.

For special weakly complete Hopf algebras $A$ which are close enough to
$\K[G]$, 
these insights allow us to formulate a concrete identifiction of the
dual $A'$ of the group Hopf-algebra $\K$ over the ground fields
$\R$ or $\C$, where again we consider the pair 
$({\mathcal V}, {\mathcal W})$ of dual categories $\mathcal V$ of 
$\K$-vector spaces and $\mathcal W$ of weakly complete $\K$-vector
spaces. 

For an arbitrary topological group $G$, let
 again $\K[G]'$ denote the topological dual  of the
 weakly complete group Hopf-algebra $\K[G]$ 
over the field $\K\in\{\R,\C\}$. 
Then  the following 
identification of the topological dual of $\K[G]$
in the category ${\mathcal V}$ provides the background of
the the topological dual of any weakly complete Hopf algebra
generated by its grouplike elements:

\begin{Theorem} \label{th:dual} 
 {\rm(a)} For an arbitrary topological group $G$, the function 
$$\omega\mapsto \omega\circ \eta_G: \K[G]' \to \RR(G,\K)$$
is an isomorphism of $\K$-vector spaces.

{\rm(b)} If $A$ is a weakly complete Hopf algebra satisfying
$\SS(A)=A$, and if $G$ is the group $\Gamma(A)$ of grouplike 
elements of $A$, then $\tau_A\colon A' \to \RR(G,\K)$ implements an
isomorphism onto a Hopf subalgebra of the Hopf algebra $\RR(G,\K)$
dual to the group Hopf algebra $\K[G]$.
\end{Theorem}

In the sense of (a) in the preceding theorem, $\K[G]'$ may be 
identified with the  vector space 
$\RR(G,\K)$ of continuous functions
$f\colon G\to \K$ for which  that the 
translates ${}_gf$ (where ${}_gf(h)=f(gh)$) 
(equivalently, the translates $f_g$ (where $f_g(h)=f(gh)$) span a finite 
dimensional $\K$-vector space.
The vector space $\RR(G,\K)$ is familiar in the literature as the
vector space of representative functions on $G$, where it is 
most frequently formulated for compact groups $G$ and where it
is also considered as a Hopf-algebra. In that case, the isomorphism
of Theorem \ref{th:dual} is also an isomorphism of Hopf algebras. 
We are choosing here the covariant group algebra $\K[G]$
to be at the center of attention and obtain $\RR(G,\K)$
via vector space duality from $\K[G]$. Conversely, if one asks for
a ``concrete'' description of $\K[G]$, then the answer may now be
that, in terms of topological vector spaces, as a topological vector 
space, $\K[G]$ is the {\color{\red}algebraic} dual {\color{\red} 
(consisting of all linear forms)} of
the (abstract) vector subspace $\RR(G,\K)$  of the 
vector space $C(G,\K)$ of continuous functions $G\to \K$.
If $G$ is a compact group, $C(G,\K)$ is a familiar
Banach space.

\bsk
At this point one realizes that in case of a finite group $G$ 
the $\K$-vector spaces $\K[G]$ and $\RR(G,\K)$ are finite dimensional
of the same dimension and are, therefore isomorphic. In particular,
we have a class of special examples of $\K$-Hopf algebras:

\begin{Example} \label{e:fin-dim} Let $G$ be a finite group. Then
$A=\RR(G,\K)=\K^G$ is a finite dimensional real Hopf algebra which
may be considered as the dual of the group algebra $\K[G]$. 
A grouplike element $f\in\Gamma(A)$ then is a  character 
$f\colon \K[G]\to \K$. The group $G$ may be considered as 
the subgroup $\Gamma(\K[G])$, and then $f|G\colon G\to \K^{-1}$
is a morphism of groups.  Conversely, every group morphism
$G\to \K^{-1}$ yields a grouplike element of $A$. In view of the
finiteness of $G$ we have
$$ \Gamma\(\RR(G,\K)\)\cong \Hom(G,\Ss^1),\quad \Ss^1=\{z\in\C:|z|=1\}.$$
If $G$ is a finite simple group, then  the $|G|$-dimensional
Hopf algebra $A=\RR(G,\K)$ has no nontrivial
grouplike elements, that is, $\Gamma(A)=\{1\}$ and $\SS(A)=\K\.1$.
 The smallest example of this kind arises from  the group $G=A_5$, the
group of 120 even permutations of a set of 5 elements.
\end{Example}

\section{$\C[G]$ as an Involutive Weakly Complete  Algebra} 

A complex unital algebra $A$ is called {\it involutive} if 
{\color{\red}it is endowed with} 
a real vector space automorphism $a\mapsto a^*$ such that for
all $a, b\in A$ with have  $a^{**}=a$, ${\color{\red}(c\.a)^*}
=\overline c\.a^*$,
and $(ab)^*=b^*a^*$. A complex weakly complete unital topological
algebra whose underlying algebra is involutive 
with respect to an involution $^*$ is a weakly complete unital
 involutive algebra, if $a\mapsto a^*$ is an isomorphism of the underlying
weakly complete topological vector space.

\hskip-10pt Every C$^*$-Algebra is an involutive algebra;
a simple example is $\Hom(\C^n,\C^n)$, i.e, an algebra isomorphic
to the algebra of $n\times n$ complex matrices with the conjugate
transpose as involution. An element $a\in A$ is {\it unitary} if
$aa^*=a^*a=1$, that is if $a^*=a^{-1}$, and it is called 
{\it hermitian} if $a=a^*$. If $A$ is an involutive topological 
unital algebra with an exponential function $\exp$, then for each
hermitian element $h$ the element {\color{\red}$\exp(i\.h)$} is unitary.
 
Note that in $A=\C$ the unitary elements are the elements on 
the unit circle $\SS^1$ and the hermitian elements are the ones
on the real line $\R$. 
The theory of the function $h\mapsto{\color{\red}\exp(ih)}$ is called {\it
trigonometry}.   

\msk

The involutive algebras form
a category with respect to morphisms $f\colon A\to B$ satisfying
$f(a^*)=f(a)^*$. 

\msk

\begin{Lemma} \label{l:generalnonsense}
If $A$ is a complex weakly complete unital topological coalgebra,
then the complex vector space dual $A'$ is
 an involutive unital algebra for the
involution $f\mapsto f^*$ defined by $f^*(a)=\overline{f(a^*)}$. 
\end{Lemma}

The proof is a straightforward exercise.
\bsk

Assume now that $G$ is a topological group.
 We shall  introduce an involution on $\C[G]$ making it an involutive
algebra.

\msk

(a) For  every complex vector space $V$ we introduce a complex vector
space $\til V$ by endowing the real vector space underlying $V$ with a
complex scalar multiplication $\bullet$ defined by 
$c{\bullet}v=\overline c\.v$. 

(b) The composition 
$$ |\C[G]|\mapright{\sigma_G}|\C[G]|\mapright{\id_{\C[G]}}
   |\C[G]|\til{\phantom m}$$
of the Hopf algebra symmetry $\sigma_G: |\C[G]|\to |\C[G]|$ and
the identity map $\C[G]\to \C[G]\til{\phantom m}$ 
yields an involution
$$^*\colon|\C[G]{\color{\red}|}\to|\C[G]|,\quad \mbox{ \it such that }
(\forall c\in\C,\ a\in\C[G])\, a^*=\sigma(a)\mbox { \it and }(c\.a)^*=
   \overline c\.a^*.$$
Moreover, by definition of $\sigma_G$, we have
$\eta_G(g)^*=\eta_G(g^{-1})$ for all $g\in  G$.

\msk

The proof of the following remarks may again be safely left as an exercise.

\begin{Lemma} \label{l:clear} For each topological group $G$,
the complex algebra $\C[G]$ is an involutive algebra with respect
to $^*$, and the comultiplication $c\colon \C[G]\to\C[G]\otimes\C[G]$
and the coidentity $k\colon\C[G]\to\C$ are morphisms of involutive 
algebras. All elements in 
$\eta_G(G)\subseteq \C[G]$ are unitary grouplike elements.
\end{Lemma}

From Lemma \ref{l:generalnonsense} we derive directly the following
observation:

\begin{Lemma} \label{l:alsoclear} The dual $\RR(G,\C)$ of the
weakly complete involutive algebra $\C[G]$ is an involutive
algebra.
\end{Lemma}

\begin{Definition} \label{d:sixfour2} Let $R$ be a unital algebra over $\C$. Then
 the set of hermitian characters
of $R$ is denoted by $\Spec_h(R)$ is called the {\it hermitian spectrum}
of $R$.
\end{Definition}

For a {\it compact} group $G$, however, 
the commutative unital algebra $\RR(G,\K)$
is a {\it dense} subalgebra of the commutative unital Banach $\K$--algebra
$C(G,\K)$. (See e.g.\ \cite{book3}, Theorem 3.7.)

We note that for every $g\in G$, the {\it point evaluation} 
$f\mapsto f(g): \RR(G,\C)\to \C$ belongs to $\Spec_h(\RR(G,\C))$.
In this particular
situation, the literature provides the following result in which we
consider $\C[G]$ as the dual $\RR(G,\C)^*$ of the involutive
algebra $\RR(G,\C)$.

\msk

\begin{Theorem} \label{th:classic} For a compact group $G$, the 
  hermitian spectrum 

\centerline{$\Spec_h(\RR(G,\C))\subseteq \C[G]$}

\nin
is precisely the set of point evaluations.
\end{Theorem}

\begin{Proof} See e.g. \cite{hofbig}, Proposition 12.26  
or \cite{hoc}, p.~28. Cf.\ also \cite{dix}, no. {\bf1.}3.7.\ on 
p.~7.
\end{Proof} 

In the case $\K=\R$, any morphism $f\colon \RR(G,\R)\to \R$ of real algebras,
that is, any (real) character of $\RR(G,\R)$ extends uniquely to 
a (complex) character $\til f\colon \C\otimes_\R \RR(G,\R)\to \C$, where
$\C\otimes_{\color{\red}\R}\RR(G,\R)\cong\RR(G,\C)$, such that 
{\color{\red}$\tilbar f= \bartil f$}. 
Trivially,  $f$ is a point evaluation of $\RR(G,\R)$, if and only if  $\til f$
is a point evaluation of $\RR(G,\C)$, and $f$ is continuous if and only if
$\til f$ is continuous. Hence Theorem \ref{th:classic} implies
the analogous result over $\R$:

\begin{Corollary} \label{c:classic} For a compact group $G$, the 
set $\Spec(\RR(G,\R))\subseteq \RR(G,\R)^*$ of
characters of $\RR(G,\R)$ is precisely the set of point evaluations
$f\mapsto f(g)$, $G\in G$.
\end{Corollary} 

{\color{\red} 
Recall that $\RR(G,\R)^*$ and $\R[G]$ are natually isomorphic by
duality since $\R[G]'$ and $\RR(G,\R)$ are naturally isomorphic
by Theorem \ref{th:dual}.

After Theorem \ref{th:gr-alg-comp-gr} we may identify a compact
group $G$ with its isomorphic image via $\eta_G$ in $\R[G]$.
From Remark \ref{r:grlike} we know $G\subseteq \Gamma(\R[G])$
(see Definition \ref{d:grlike}). Now finally we have
the full information on one half of the adjunction of
$\HH$ and $\Gamma$ of Theorem \ref{p:grounding} in the case
of {\it compact} groups $G$:

\begin{Theorem} \label{th:comp-group-like} {\rm(Main Theorem for Compact Groups A)}
\quad  For every compact group $G$
the natural morphism of topological groups $\eta_G\colon G\to \Gamma\(\HH(G)\)
=\Gamma(\R[G])$ is an isomorphism of compact groups.
\end{Theorem}

In other words, if for a compact group $G$ we identify $G$ with 
a subgroup of $(\R[G],\.)$ for its weakly complete group
algebra $\R[G]$, then $G=\Gamma(\R[G])$, i.e. every grouplike
element of the Hopf algebra $\R[G]$ is a group element of $G$.

\begin{Proof} We apply Proposition \ref{p:grouplike} with 
$R=\RR(G,\R)$ and $A=\R[G]$. Then the set of grouplike elements
of $\R[G]$ is $\Spec(\RR(G,\R))$. Let $\ev\colon G\to\Spec(\RR(G,\R))$
denote the function given by $\ev(g)(f)=f(g)$ for all $g\in G$,
$f\in\RR(G,\R)$. Then $\ev(G)\subseteq \Spec(\RR(G,\R))$ by
Remark \ref{r:grlike} and Proposition \ref{p:grouplike}.
However, Theorem \ref{th:classic} and Corollary \ref{c:classic}
show that in fact equality holds, and that is the assertion of
the theorem.
\end{Proof}

It appears that a direct proof of the assertion that everery
grouplike element is a group element in $\R[G]$ would be
difficult to come by even in the case of an infinite compact group $G$.
(If $G$ happens to be finite then the proof is
elementary linear algebra.). 

As it stands, after Theorem \ref{th:comp-group-like}, the
category theoretical background on the adjunction of the functors
$\Gamma$ and $\HH$ expressed in Corollary \ref{c:alternate2} (1)
and (2) gives some additional insight for compact $G$, respectively,
compact $\Gamma(A)$.
Indeed Corollary \ref{c:alternate2}(1) yields at once:

\begin{Corollary} \label{c:nochwas} Let $A$ be a weakly complete real 
Hopf algebra and assume that the group $G\defi\Gamma(A)$ of its grouplike
elements is compact. Then the natural morphism of topological groups
$\Gamma(\epsilon_A)\colon \Gamma\(\HH(G)\)\to G$ is an isomorphism.
\end{Corollary}                                                                                                
 }

\bsk

We recall that every compact group $G$ is a pro-Lie group with a 
Lie algebra $\L(G)$ and an exponential function $\exp_G\colon\L(G)\to G$
according to \cite{book3}.
\msk

\begin{Theorem} \label{th:comp-prolie} {\color{\red}
{\rm (Main Theorem for Compact 
Groups B)}} \quad Consider the compact group $G$
as a subgroup of the multiplicative monoid $(\R[G],\cdot)$ of its 
weakly complete group Hopf algebra. Then {\color{\red} $G=\Gamma(\R[G])$, that is, 
$G$ is the set of grouplike elements of the Hopf algebra $\R[G]$,
 while  $\L(G)\cong\Pi(\R[G])$, that is, the Lie algebra of $G$
is the set of primitive elements of $\R[G]$.
Hence the Lie algebra of $G$ is isomorphic
to the  subalgebra $\Pi(\R[G])$ of the Lie algebra 
$\R[G]_{\rm Lie}$.} Moreover, the global exponential function 
$\exp_{{\color{\red}\R}[G]}:\R[G]_{\rm Lie} \to \R[G]^{-1}$
of $\R[G]$ restricts{\color{\red}, up to isomorphy,} 
to the exponential function $\exp_G\colon\L(G)\to G$
of $G$.
\end{Theorem}

\msk

We have now seen that a compact group and its Lie theory are completely
incorporated into its weakly complete group Hopf algebra. 
Under such circumstances it is natural to ask whether a similar assertion
could be made for the ample theory of Radon measures on a compact
group including its
Haar measure.
We shall address this question and give a largely affirmative answer
 in the next section.

\msk

 What we have to address  at this point is the question, whether for
a weakly complete real Hopf algebra $A$ in which the group
$\Gamma(A)$ is compact and algebraically and topologically generates $A$,
 the natural morphism
$\epsilon_A\colon \HH\(\Gamma(A)\)=\R[\Gamma(A)]\to A$
is in fact an isomorphism.

\begingroup
\color{\red} For the investigation of this question we need some preparation.
Assume that $G$ is a compact group and 
$\RR(G,R)\subseteq C(G,\R)$ is the Hopf algebra of all
functions $f\in C(G,\R)$ whose translates span a finite dimensional 
vector subspace.

We now  let $M$  be a Hopf subalgebra and a
$G$-submodule of $\RR(G,\R)$. Recall that $\Spec M$ denotes the set of all
algebra morphisms $M\to\R$. Trivially we have  a morphism
 $\omega \mapsto \omega|M:\Spec \RR(G,\R) \to \Spec(M)$.
 From Corollary \ref{c:classic} we know that 
$G \cong \Spec \RR(G,\R)$ via point evaluation. 
In the case that $M=A'$ and  $G=\Gamma(A)$ such that 
$\Gamma(\epsilon_A)\colon \Gamma(\R[G])\to G$  an isomorphism
as in Corollary \ref{c:nochwas}
 we know that
$$g\mapsto(f\mapsto f(g)):G\to \Spec(M)\quad\mbox{is an isomorphism.}
\leqno(E)$$
From \cite{book3}, Definition 1.20, p.~13 we recall that 
$M\subseteq C(G,\R)$ is said to {\it separate points} if
for two points $g_1\ne g_2$ in $G$ there is an $f\in M$ such that
$f(g_1)\ne f(g_2)$. In other words, different points in $G$ can 
be distinguished by different point evaluations of functions from
$M$. So condition (E) secures that the functions of $M$ separate
the points of $G$. 
This has the following consequence:

\begin{Lemma} \label{l:stone-wei} If the unital subalgebra $M$ of $\RR(G,\R)$
 satisfies $(E)$, then 
$M$ is dense in $C(G,\R)$ with respect to the sup norm 
 and  is dense in $L^2(G,\R)$ with respect to the $L^2$-norm.
\end{Lemma}

\begin{Proof} Since $M$ is a unital subalgebra of $\RR(G,\R)$, 
it  contains the scalar multiples of
 the constant functions of value 1, that is, $M$ contains all
the constant functions.
Moreover,  by Hypothesis (E), 
the algebra $M\subseteq C(G,\R)$
separates the points of $G$. Therefore the Stone-Weierstra{\ss} Theorem
applies and shows that $M$ is dense in $C(G,\R)$ in the sup norm topology
of $C(G,\R)$.  
Since $L^2(G,\R)$ is the $L^2$-norm completion of
$C(G,\R)$ and $M$ is uniformly dense in $C(G,R)$ 
 it follows
that $M$ is dense in $L^2(G,\R)$ in the $L^2$-norm.
(Cf.\  e.g.\ \cite {book3} Theorem 3.7 and its proof.)
\end{Proof}

\begin{Lemma} \label{l:crucial}  If a $G$-submodule $M$ of $\RR(G,E)$ 
is $L^2$-dense in $L^2(G,\R)$ then it agrees with $\RR(G,R)$.
\end{Lemma}

\begin{Proof} Let $\hat G$ denote the set of isomorphy classes 
of irreducible $G$-modules. By the Fine Structure Theorem of $\RR(G,\R)$
 (see \cite{book3}, Theorem 3.28) 
 $\RR(G,R)=\sum_{\epsilon\in\hat G}\RR(G,\R)_\epsilon$
where $\sum$ denotes the algebraic direct sum of
(finite dimensional) vector subspaces and where
 $\RR(G,\R)_\epsilon$ is a finite direct sum of simple
modules for
each $\epsilon\in \hat G$. In particular, each $\RR(G,\R)_\epsilon$ 
is finite dimensional.
Further $L^2(G,\R)=\bigoplus_{\epsilon\in\hat G} \RR(G,\R)_\epsilon$ where 
$\bigoplus$ denotes the Hilbert space direct sum.

The submodule $M$ of $\RR(G,\R)$ adjusts to the canonical decomposition of
$\RR(G,\R)$ since $M_\epsilon$ is necessarily a submodule of 
$\RR(G,\R)_\epsilon$. Hence
$M=\sum_{\epsilon\in\hat G}M_\epsilon$ and the 
$L_2$-closure of $M$ in $L^2(G,\R)$ is
 the Hilbert space sum $\bigoplus_{\epsilon\in\hat G} M_\epsilon$.

\msk

By way of contradiction suppose now that that $M\ne \RR(G,\R)$.
Then there is an $\epsilon'\in\hat G$ such that 
$M_{\epsilon'}\ne \RR(G,\R)_{\epsilon'}$. Since
$\RR(G,\R)_{\epsilon'}$ is finite dimensional,
$$\bigoplus_{\epsilon\in\hat G}M_\epsilon=M_{\epsilon'}\oplus
\bigoplus_{\epsilon\ne\epsilon'}M_\epsilon$$
is properly smaller than 
the Hilbert space sum 
$$\bigoplus_{\epsilon\in\hat G}\RR(G,\R)_\epsilon
=\RR(G,\R)_{\epsilon'}\oplus\bigoplus_{\epsilon\ne\epsilon'}\RR(G,R)_\epsilon
=L^2(G,\R),$$ contradicting the hypothesis that $M$ is $L^2$-dense
in $L^2(G,\R)$.
This contradiction proves the lemma.
\end{Proof}

Now we are ready for the third main result on compact groups
in the present context: the statement which parallels Theorem 
\ref{th:comp-group-like}.

\begin{Theorem} \label{th:comp-prolie2} {\rm (Main Theorem on Compact
Groups C)}\quad  Let $A$ be a weakly complete real Hopf Algebra satisfying
the following two conditions:
\begin{enumerate}[\rm(i)]
\item The subgroup 
$\Gamma(A)$ of grouplike elements of $A$ is compact,
\item  $\Gamma(A)$ generates $A$ 
algebraically and topologically, that is, $\SS(A)=A$.
\end{enumerate}
 Then 
$$\epsilon_A\colon \R[\Gamma(A)]=\HH\(\Gamma(A)\)\to A$$ is a natural
isomorphism.
\end{Theorem} 

\begin{Proof} We set $G=\Gamma(A)$. By Corollary \ref{c:epsilon-epic} the morphism
$\epsilon_A\colon \R[G]\to A$ is a quotien homomorphism of weakly complete
Hopf algebras which by Corollary \ref{c:nochwas} 
induces an isomorphism $\Gamma(\epsilon_A)\colon \Gamma(\R[G])\to G$.
So $\Gamma(\R[G])$ is identified with $G$ if we consider $G$ as included in $\R[G]$
according to Theorem \ref{th:comp-group-like}.

By the Duality between real Hopf algebras in $\mathcal V$ and 
weakly complete real Hopf algebras in $\mathcal W$, the dual morphism
$\epsilon_A'\colon A'\to \R[G]'$ is injective. Theorem \ref{th:dual} 
then gives us an inclusion  $A'\subseteq\RR(G,\R)$ of real Hopf algebras 
as well as of $G$-modules 
such that the natural map $\Spec(A')\to\Spec\(\RR(G,\R)\)$ is
the identity.

That is, Condition (E) above (preceding Lemma \ref{l:stone-wei})
holds and Lemmas \ref{l:stone-wei} and \ref{l:crucial} apply.
Therefore $A'=\RR(G,\R)$. This in turn shows that $\epsilon_A$
is an isomorphism. 
\end{Proof} 

For a concise formulation of the consequences let us use
the following notation:

\begin{Definition} A real weakly complete Hopf algebra $A$ will be called 
{\it compactlike} if it is  such that
the subgroup $\Gamma(A)$ of grouplike elements is compact
and $\SS(A)=\overline{\Span}\(\Gamma(A)\)=A$.
\end{Definition} 

\begin{Theorem} {\rm (The Equivalence Theorem of the Category
of Compact Groups and compactlike Hopf algebras)} The categories of compact groups and of
weakly complete compactlike Hopf algebras are equivalent.
\end{Theorem}

\begin{Proof} This follows immediately from Theorems
\ref{th:comp-group-like} and \ref{th:comp-prolie2}.
\end{Proof}
 
\begin{Corollary} \label{c:tannaka} {\rm(Tannaka Duality)}\quad
The category of compact groups
is dual to the full category of real abstract Hopf algebras of the
form $\RR(G,\R)$ with a compact group $G$.
\end{Corollary}

\endgroup

\section{The  Radon Measures within the
           Group Algebra of a Compact Group}

We shall invoke measure theory in the form  pioneered for arbitrarly locally
compact groups in \cite{boui}. For a {\it compact} group $G$ it is 
less technical and adapts reasonably to the formalism of its real group
algebra $\R[G]$. This discussion will help us to understand the
power ot the group algebras $\R[G]$ for a compact group.

\ssk

Indeed any compact Hausdorff topological group provides us with a real Banach algebra
$C(G,\R)$ endowed with the sup-norm which makes it into  a Hopf-algebra in the categoy
of Banach spaces. Accordingly,
its topological dual $C(G,\R)'$ yields the Banach algebra and
indeed Banach Hopf algebra  $M(G,\R)$ (see e.g.\ \cite{hofbig}). 
Its elements  $\mu$ are 
the so called  {\it Radon measures} on $G$ and the general source books
of this orientation of  measure and probability theory is Bourbaki's
book \cite{boui} and for the foundations of harmonic analysis
the book of Hewitt and Ross \cite{hewross}.
For a measure theory in the context of compact groups see 
also \cite{book3}, Appendix 5: ``Measures on Compact Groups''.

\msk

So let $W$ be a weakly complete real
vector space. Then $W$ may be identified with ${W'}^*$ (see Theorem
\ref{th:vect}). For $F\in C(G,W)$ and $\mu\in M(G,\R)$
we obtain a unique element  $\int_G F\, d\mu\in W$
such that  we have
$$(\forall \omega\in W')\quad \left\<\omega,\int_GF\,d\mu\right\>
=\int_{g\in G}\<\omega,F(g)\>d\mu(g).\leqno(*)$$
(See \cite{boui}, Chap. III, \S 3, {\bf n}$^{{\rm o} 1}$, D\'efinition 1.)
Let $\supp(\mu)$ denote the support of $\mu$. (See \cite{boui}, 
Chap. III, \S2, {\bf n}$^{{\rm o} 2}$, D\'efinition 1.)

\msk

\begin{Lemma}\label{l:lin} Let $T\colon W_1\to W_2$ be a morphism of
weakly complete vector spaces, $G$ a compact Hausdorff space
and  $\mu$  a measure on $G$. 
If $F\in C(G,W_1)$, then $T(\int_G F\,d\mu)=\int_G (T\circ F) d\mu$.
\end{Lemma}
\noindent(See e.g.\ \cite{boui}, 
Chap. III, \S 3, {\bf n}$^{{\rm o} 2}$, Proposition 2.)

In \cite{boui} it is shown that   the
vector space $M(G,\R)$ is also a complete lattice w.r.t.\ a 
natural partial order 
(see \cite{boui}, Chap. III, \S 1, {\bf n}$^{\rm o}$ 6)
so that each $\mu\in M(G)$ is uniquely of the form $\mu=\mu^+-\mu^-$
for the two positive measures $\mu^+=\mu\vee 0$ and $\mu^-=-\mu\vee 0$.
 One defines $|\mu|=\mu^++\mu^-$. If $M^+(G)$ denotes the cone of all
positive measures, we have $M(G)=M^+(G)-M^+(G)$ (\cite{boui}, Chap. III,
\S 1, {\bf n}$^{\rm o}$ 5, Th\'eor\`eme 2). 
Moreover, $\|\mu\|=|\mu|(1)=\int d|\mu|$.
A measure is called a {\it probability measure} if it is positive 
and $\mu(1)=1$.
We write $P(G)$ for the set of all probability measures on $G$ and
we note $M^+(G)=\R_+\.P(G)$ where $\R_+=[0,\infty[\ \subseteq\R$.
We denote by $M_p(G)$  the vector space $M(G,\R)$ with
the topology of pointwise convergence and
 recall that $P(G)$ has the structure of a compact submonoid of
$M_p(G)^\times$; some aspects are discussed in $\cite{book3}$,
Appendix 5.
On $M^+(G)$ the topologies of $M_p(G)$ and the compact open topology 
of $M(G,\R)$ agree
(\cite{boui}, Chap. III, \S 1, {\bf n}$^{\rm o}$ 10, Proposition 18),
Also $M^+_p(G)$ is a locally compact convex pointed {\it cone} 
with the closed convex hull $P(G)$ of the set of point 
measures as {\it basis}.           
We also recall, that any positive linear form on $C(G,\R)$ is
in $M^+(G)$ (i.e., is continuous) (see \cite{boui},
Chap. III, \S 1, {\bf n}$^{\rm o}$ 5, Th\'eor\`eme 1).

\bsk

\subsection{Measures and  Group Algebras}

Now we allow this machinery and 
Theorems \ref{th:hopf-prolie} and \ref{th:gr-alg-comp-gr}
to come together to elucidate the structure of $\R[G]$ for compact groups $G$. 

\bsk
We let $G$ be a compact group.   By Theorem \ref{th:gr-alg-comp-gr}
it is no loss of generality
to assume that $G$ is a compact subgroup of  $\R[G]^{-1}$, where $\R[G]$ is 
the weakly complete group Hopf algebra of $G$ and 
$\eta_G\colon G\to \R[G]^{-1}$
is the inclusion morphism.
By Theorem \ref{th:dual} there is an isomorphism
$\omega\mapsto f_\omega\colon \R[G]'\to \RR(G,\R)$ such that 
$$(\forall \omega\in\R[G]',\, g\in G)\, \<\omega,\eta_G(g)\>=f_\omega(g),$$
and, in the reverse direction,   the function 
$\omega\mapsto \omega|G: \R[G]'\to C(G,\R)$
induces an isomorphism of vector spaces 
$\R[G]'\to{\mathcal R}(G,\R)$.

\msk

Therefore, in the spirit of relation $(*)$, we 
are led to the following definition

\msk

\begin{Definition} \label{d:maindef}
Let $G$ be a compact group. Then
each $\mu\in M(G,\R)$ gives rise to an element 
$$\rho_G(\mu)\defi \int_G\eta_G\,d\mu\in \R[G]$$
such that for all $\omega\in\R[G]'$ we have
$$\<\omega,\rho_G(\mu)\>=\int_{g\in G} \<\omega,\eta_G\>\,d\mu(g)
=\int_{g\in G}f_\omega(g)\,d\mu(g)=\mu(f_\omega).\leqno(**)$$
Therefore we have a morphism of  vector spaces
$$\rho_G\colon M(G,\R)\to\R[G].$$

We let $\tau_{\RR(G,\R)}$ denote the weakest topology making the
functions $\mu\mapsto \mu(f):M(G,\R)\to\R$ {\color{\red}continuous}
for all $f\in\RR(G,\R)$
\end{Definition}

On any compact subspace of $M_p(G)$ such as $P(G)$ the topology
$\tau_{\RR(G,\R)}$ agrees with the topology of $M_p(G)$.
\msk

\begin{Lemma} \label{l:injective} The morphism 
$\rho_G$ is injective and has dense
image.
\end{Lemma}
\begin{Proof}
 We observe $\mu\in\ker\rho_G$ if
for all $f\in\RR(G,\R)$ we have $\int_{g\in G}f(g)\,d\mu(g)=0$.
Since $\mu$ is continuous on $C(G,\R)$ in the norm topology and 
$\RR(G,\R)$ is dense in $C(G,\R)$ by the Theorem of Peter and
Weyl (see e.g.\ \cite{book3}, Theorem 3.7), it follows that
$\mu=0$. So $\rho_G$ is injective.

 If $\mu=\delta_x$ is a measure with support $\{x\}$ for
some $x\in G$, then $\rho_G(\mu)=\int_G\eta_G d\delta_x=x$.
Thus $G\subseteq\rho_G(M(G))$. Since $\R[G]$ is the closed linear
span of $G$ by Proposition \ref{p:generation}, 
it follows that $\rho_G$ has a dense image. 
\end{Proof}

We note that in some sense $\rho_G$ is 
dual to the inclusion morphism of  vector spaces
$\sigma_G\colon\RR(G,\R)\to C(G,\R)$.

Returning to (**) in Definition \ref{d:maindef}, for a compact
group $G$, we observe

\begin{Lemma} \label{l:topologies} The morphism
$$\rho_G:(M(G,\R),\tau_{\RR(G,\R)}) \to \R[G]$$ 
is a topological embedding.
\end{Lemma}

\msk

If $\mu$ is a probabililty measure,
then the element $\rho_G(\mu)=\int_G\eta_G\,d\mu$ is
contained in the compact closed convex hull
$\overline{\conv}(G)\subseteq \R[G]$. Intuitively, 
$\int_G \eta_G d\mu\in\overline{\conv}(G)$ is the center of gravity of
the ``mass'' distribution  $\mu$ contained in $G\subseteq\R[G]$. In particular,
if  $\gamma\in M(G,\R)$ denotes normalized Haar measure on $G$, then 
$$\rho_G(\gamma)=\int_G\eta_G\,d\gamma=\int_{g\in G} g\, dg$$
 is the center of gravity
of $G$ itself with respect to Haar measure.

\msk 

We note that in the weakly complete vector space $\R[G]$ the
closed convex hull 
$$ B(G)\defi \overline{\conv}(G) \subseteq \R[G]$$
is compact. (See e.g.\ \cite{book3}, Exercise E3.13.)

\begin{Lemma}\label{l:conv} The restriction 
$\rho_G|P(G): P(G)\to B(G)$
is an  affine homeomorphism.
\end{Lemma} 
\begin{Proof} (i) Affinity is clear and injectivivity
we know from Lemma \ref{l:injective}{\color{\red}.} 

(ii) Since $P(G)$ is compact in the weak topology
and $\rho_G$ is injective and continuous, $\rho_G|P(G)$
is a homeomorphism onto its image. But $G\subseteq \rho_G(P(G))$,
 and $B(G)$ is the closed convex hull of $G$ in $\R[G]$,
it follows that $B(G)\subseteq\rho_G(p(G))$.
\end{Proof}

\msk

If $k\colon \R[G] \to \R$ is the augmentation map 
(i.e., the coidentity morphism),
then $k(G)=\{1\}$ and so $k(B(G))=\{1\}$ as well. 
From $GG\subseteq G$ we deduce that
$\conv(G)\conv(G)\subseteq\conv(G)$ and from there, 
by the continuity of the
multiplication in $\R[G]$ and $1\in G\subseteq B(G)$, it follows
that $B(G)$ is a compact submonoid of $\R[G]^\times$ 
contained in the submonoid $k^{-1}(1)$.

\msk

Then the cone $\R_+[G]\defi\R_+\.B(G)$,
due to the compactness of $B(G)$,  is a locally compact submonoid
as well. The set 

\centerline{$k^{-1}(1)\cap \R_+[G]=\{{\color{\red}x\in \R_+[G]:k(x)}=1\}
=B(G)$}

\nin is a compact basis of the cone $\R_+[G]$. 

\msk

\begin{Corollary} \label{c:conv} The function
$\rho_{{\color{\red}G}}|M^+(G):M^+(G)\to \R_+[G]$ is an isomorphism of 
convex cones and 
$\rho_G(M(G))=\R_+[G]-\R_+[G]$
\end{Corollary}

\begin{Proof} Since $M^+(G)=\R_+\.P(G)$ and $\R_+[G]=\R_+\.B(G)$,
Lemma \ref{l:conv} shows that $\rho_G|M^+(G)$ is an affine 
homeomorphism. Since $M(G)=M^+(G)-M^+(G)$, the corollary follows.
\end{Proof}

Among other things  this means that every element of $\R_+[G]-\R_+[G]$
is a an integral $\int_G \eta\,d\mu$ in $\R[G]$ for some Radon 
measure $\mu\in M(G)$ on $G$.

\msk

In order to summarize our findings 
 we  firstly list the required conventions:

\nin
Let $G$ be   a compact group viewed as
a subgroup of the group $\R[G]^{-1}$ of units of the 
weakly complete group algebra  $\R[G]$. Let
$B(G)=\overline{\conv}(G)$ denote the closed convex hull of $G$ in $\R[G]$
and define $\R_+[G]=\R_+\.B(G)$. Let $k\colon \R[G]\to \R$ denote the
augmentation morphism and $I=\ker k$ the augmentation ideal. 
We let $\eta_G\colon G\to\R[G]$ denote the inclusion map and
consider $\rho_G\colon M(G,\R)\to\R[G]$ with 
$\rho_G(\mu)=\int_G\eta_G\,d\mu$.

\begin{Theorem} \label{th:convexity} For a compact group $G$
we have the following conclusions:
\begin{enumerate}[\rm (a)]
\item  $B(G)\supseteq G$ is a compact submonoid of 
       $1+I\subseteq {\color{\red}(\R[G],\cdot)}$
       with Haar measure $\gamma$ of $G$ as zero element.

\item  $\R_+[G]$ is a locally compact pointed cone with basis 
       $B(G)$, and is
       a  submonoid of ${\color{\red}(\R[G],\cdot)}$. 

\item  The function $\rho_G\colon (M(G),\tau_{\RR(G,\R)})\to\R[G]$ 
       is an injective morphism of topological vector spaces
       {\color{\red} with dense image $\R_+[G]-\R_+[G]$. It
       induces a homeomorphism onto its image.}

\item  The function $\rho_G|M^+(G)\colon M_p^+(G)\to \R_+[G]$ is an 
       affine homeomorphism from
       the  locally compact convex cone  
        of positive Radon measures on $G$
       onto $\R_+[G]\supseteq B(G)\supseteq G$.
 \end{enumerate}
\end{Theorem}

\begin{Remark} The Haar measure $\gamma$ is  mapped
by $\rho_G$ onto the center of gravity $\int_G \eta_G\,d\gamma$
of {\color{\red}$G$, $\gamma\in B(G)\subseteq 1+I$}.
\end{Remark}
\bsk
It should be noted that $\rho_G\colon M(G,\R)\to\R[G]$ is far
from surjective if $G$ is infinite: 
If we identify $\R[G]$ with $\RR(G,\R)^*$
according to Theorem \ref{th:dual}, then any element $u\in\R[G]$
representing a {\color{\red} linear form on $\RR(G,\R)$ which is
discontinuous in the norm topology induced by $C(G,\R)$}
fails to be an element of $\rho_G(M(G))$.

\msk 
Theorem \ref{th:convexity} shows that for a compact group $G$,
the weakly complete real group  algebra $\R[G]$ does not only
contain $G$ and the entire pro-Lie group theory encapsulated
in the exponential function $\exp_G\colon \L(G)\to G$ but also
the measure theory, notably, that of the monoid of probability
measures $P(G)\cong B(G)$.

\bsk

Recall the hyperplane ideal $I=\ker k$ for the augmentation
$k\colon \R[G]\to\R$.

\begin{Corollary} \label{c:split} Let $G$ be a compact group, $\R[G]$ its real
group algebra, and $\gamma\in\R[G]$ its normalized Haar measure.
 Then $J\defi \R\.\gamma$ is a one-dimensional ideal,  and
$$\R[G]=I\oplus J$$
is the ideal direct sum of $I$ and $J$. The vector subspace $J$ is a minimal
nonzero ideal.

In particular, $J\cong \R[G]/I\cong\R$ and $I\cong \R[G]/J$.
\end{Corollary}

\begin{Proof} In the multiplicative monoid 
$B(G)\subseteq \R[G]${\color{\red},}
the idempotent element $\lambda$ is a {\color{\red} zero
of the monoid $(I,\cdot)$}, that is, 
$$\lambda B(G)=B(G)\lambda =\{\lambda\}.$$
(See \cite{book3}, Corollary A5.12.) As a consequence,
$$ JB(G) = B(G)J\subseteq J.$$
The vector space $\Span B(G)$ contains $\Span G$ which is
dense in $\R[G]$ by Proposition \ref{p:generation}.
Hence $J\R[G]=\R[G]J\subseteq J$ and so $J$ is a two-sided
ideal. Since $k(B(G))=\{1\}$ by Theorem \ref{th:convexity}(a)
we know $J\not\subseteq I$, and since $I$ is a hyperplane,
$\R[G]=I\oplus J$ follows.
\end{Proof}

\msk
We note that $\R[G]/J$ is a weakly complete topological algebra containing
a copy of $G$  and indeed of $P(G)$ with Haar measure in the copy of
$P(G)$ being the zero of the algebra. It is not a group algebra nor
a Hopf algebra in general
as the example of $G=\Z(3)$ shows. 

{\color{\red}

\msk

While the group $\Gamma(\R[G])\cong G$ of grouplike elements of $\R[G]$
(and its closed convex hull $B(G)$) is contained in the affine hyperplane
$1+I$, in the light of Theorem \ref{th:comp-prolie} it is appropriate to
observe that  in the circumstances of {\rm Corollary \ref{c:split}},
the Lie algebra of primitive elements $\Pi(\R[G])\cong \L(G)$ is
contained in $I=\ker k$.

Indeed the ground field $\R$ is itself a Hopf algebra 
with the natural isomorphism $c_\R\colon \R\to \R\otimes \R$ 
satisfying $c_\R(r)=
r\.(1\otimes1)=r\otimes1=1\otimes r$. Now
the coidentity $k$ of any coalgebra $A$ is  a morphism of coalgebras
so that we have  a commutative diagram for $A\defi\R[G]$: 

\vglue-10pt

$$\begin{matrix} A&\mapright{c_A}&A\otimes A\\
          \lmapdown{k}&&\mapdown{k\otimes k}\\
                 \R&\mapright{c_\R}&\R\otimes\R,\end{matrix}\leqno(1)
$$
If $a\in A$ is primitive, then $c(a)=a\otimes1 + 1\otimes a$.
The commutativity of (1) provides $k(a)\otimes1=\alpha(k(a))
(k\otimes k)(c(a))=(k\otimes k)(a\otimes1 +1\otimes a)=
k(a)\otimes1 + 1\otimes k(a)$, yielding $1\otimes k(a)=0$,
that is $k(a)=0$ which indeed means $a\in\ker k=I$.
We note that these matters are also compatible with the Main Theorem
for Compact Groups B \ref{th:comp-prolie} insofar as,
trivially, $\exp(I)\subseteq 1+I$.
}  

\bsk

\vglue25pt

\nin
{\bf Acknowledgments.}\quad The authors thank {\sc G\'abor Luk\'acs} for
extensive and deep discussions inspiring much of the present 
investigations, which could have hardly evolved without them.

\nin They also  gratefully acknowledge
the kind assistance of {\sc Mahir Can} of Tulane University, New Orleans
in the compilation of the references, notably with \cite{boch}, 
\cite{bresar}, and \cite{good}. Lemma \ref{l:canslemma} is
rightfully called Can's Lemma.

\nin We are also very grateful to the
thorough and expeditious input of the referee who has
assisted us in eliminating numerous formal flaws in
our manuscript and who directed us to relevant
points that needed further exposition.

\end{document}